\documentclass[twocolumn]{autart}    

\usepackage{amssymb,amsmath,amsopn}
\usepackage{amsfonts}
\usepackage{amstext}
\usepackage{empheq}
\usepackage{epsfig}
\usepackage{color}
\usepackage{graphicx,colordvi}
\usepackage{epstopdf}
\usepackage{float}
\usepackage{CJK}
\usepackage{cite}

\newtheorem{lemma}{ Lemma}[section]
\newtheorem{theorem}{ Theorem}[section]
\newtheorem{definition}{ Definition}[section]

\newtheorem{remark}{ Remark}[section]
\newtheorem{example}{ Example }[section]

\newtheorem{proposition}{ Proposition}[section]          

\begin{document}

\begin{frontmatter}

\title{Decentralized Strategies for Finite Population Linear-Quadratic-Gaussian Games and Teams\thanksref{footnoteinfo}} 

\thanks[footnoteinfo]{This work was supported by the National Natural Science Foundation of China
under Grants 61633014, 61773241, 61877057, 62122043, U1701264, and Australian Research Council under Grant DP200103507. 
}

\author[sdu]{Bing-Chang Wang}\ead{bcwang@sdu.edu.cn},
\author[sdu]{Huanshui Zhang}	\ead{hszhang@sdu.edu.cn},
\author[un]{Minyue Fu}\ead{minyue.fu@newcastle.edu.au},
\author[sdu]{Yong Liang}\ead{yongliang@mail.sdu.edu.cn}

\address[sdu]{School of Control Science and Engineering,
Shandong University, Jinan 250061, P. R. China.}	

\address[un]{the School of Electrical Engineering and Computing Science, University of Newcastle, NSW 2308, Australia.  }         

\begin{keyword}                           
Mean-field game,  decentralized Nash equilibrium, finite population, non-standard FBSDE,  weighted cost             
\end{keyword}                             

\begin{abstract}                          
 This paper is concerned with a new class of mean-field games which involve a finite number of agents. 
Necessary and sufficient conditions are obtained for the existence of the decentralized open-loop Nash equilibrium in terms of non-standard forward-backward stochastic differential equations (FBSDEs). 
By solving
 the FBSDEs, we design a set of decentralized
strategies by virtue of two differential Riccati equations. Instead of the 
$\varepsilon$-Nash equilibrium in classical mean-field games, the set of decentralized strategies is
shown to be a  Nash equilibrium.  For the infinite-horizon problem, a simple condition is given for the solvability of the algebraic Riccati equation arising from consensus. Furthermore, the social optimal control problem is studied. Under a mild condition, the decentralized
social optimal control and the corresponding social cost are given.
\end{abstract}

\end{frontmatter}

\section{Introduction}

The mean-field game has drawn intensive research attention because it provides an effective theoretical scheme for analyzing the collective behavior of large population multiagent systems (MASs). This has found wide applications in various disciplines, such as economics, biology, engineering, and social science \cite{BFY13, GS13, C14, CD18, weintraub2008markov, WH15}. 
Mean-field games were initiated by two groups independently. Huang \emph{et al.} designed an $\epsilon$-Nash equilibrium for a decentralized strategy with discount costs based on a Nash certainty equivalence (NCE) approach \cite{HCM07}. Independently, Lasry and Lions introduced a mean-field game model and studied the well-posedness of coupled partial differential equation systems \cite{LL07}. 
The NCE approach can be extended to cases with long run average costs \cite{LZ08} or with Markov jump parameters \cite{WZ12}.

\subsection{Literature review}
{Depending on the state-cost setup of a mean-field game, it can be classified into linear-quadratic-Gaussian (LQG) games or more general nonlinear ones. The LQG game 
is commonly adopted in mean-field studies because of its analytical tractability and close connection to practical applications. Relevant works include \cite{HCM07, LZ08, WZ13, BSYY16, HH16, MB17}. 
In contrast, a nonlinear mean-field game 
enjoys its modeling generality  
(see e.g. \cite{HMC06, LL07, CD13, CD18}).
Besides, depending on their system hierarchy, mean-field games can be classified into homogeneous, heterogeneous, or mixed.
See \cite{H10, WZ12, BLP13} for mixed games.} 

Apart from noncooperative games, mean-field social optimal control has also drawn increasing attention recently. In a social optimum problem, all players cooperate to optimize the social cost---the sum of individual costs. Social optima are linked to a type of team decision \cite{H80} but with highly complex interactions.
 The work \cite{HCM12} 
studied social optima in mean-field LQG control, and provided an asymptotic team-optimal solution. Authors in \cite{AM15} considered team-optimal control with finite population and partial information.
 For further literature, see  \cite{HN16} for socially optimal control for major-minor systems, \cite{WZ17} for the team problem with a Markov jump parameter as common random source,
 \cite{SNM18} for dynamic 
collective choice by finding a social optimum, \cite{SY19} for stochastic dynamic teams and their mean-field limit. 

\subsection{Motivation and contribution}
Most works on mean-field games and control focused on large/infinite population MASs and obtained approximate Nash equilibria. 
{However, in many practical situations (such as oligopolistic markets),
 small population or  moderate population systems are considered.} How about the case of finite population? Which type of Nash equilibria can we obtain? 
In this paper, we investigate mean-field games and teams for a \emph{finite} number of agents. 
Different from classical mean-field games, the considered systems are not confined to large-population MASs, and the population state average is generalized to the weighted sum of agents' states. The weighted average cost
may appear in 
locality-dependent models \cite{HMC10}, deep structured teams \cite{AA20}, and graphon mean-field games \cite{CH18}.

{
For a finite-population game or team,
agents access different information sets, thus the problem is \emph{decentralized} and different from classical vector optimization with a centralized designer.
By the terminology of  Witsenhausen \cite{W68}, agents have nonclassical information structures.
Due to essential difficulty arising from different information structures and complex strategic interactions, 
there is no systematic approaches to address the decentralized optimal control problem.
For mean field games, asymptotic decentralized Nash equilibria are obtained by state aggregation technique, where the difficulty is circumvented ingeniously by taking mean field approximation.}

{In this paper, we construct decentralized strategies by solving nonstandard forward-backward stochastic differential equations (FBSDEs). For finite-population LQG games, we first derive necessary and sufficient conditions for the existence of the open-loop Nash equilibrium in terms of 
FBSDEs by variational analysis. 
Due to accessible information restriction, the decentralized Nash equilibrium is given by the conditional expectation of costates (solutions to adjoint equations). This leads to a set of nonstandard FBSDEs. The key step of strategy design is to solve these equations.}
 We first construct auxiliary FBSDEs and establish the equivalent relationship for two set of FBSDEs. By decoupling auxiliary FBSDEs, we design a set of decentralized
strategies in terms of two differential Riccati equations.
Instead of the 
$\varepsilon$-Nash equilibrium, the set of decentralized
strategies is shown be an (exact) Nash equilibrium.
 For the infinite-horizon problem, we give some criterion for solvability of the algebraic Riccati equation arising from consensus problems.
Particularly, for the model of multiple integrators, agents can reach the mean-square consensus. {Furthermore, the finite-population team problem is also studied, where the social cost is a weighted sum of individual costs. By virtue of symmetric Riccati equations, we give the decentralized social optimal control and the corresponding social cost. In addition, the result is extended to the case with state coupling, where the strategic interaction is more complex and the auxiliary FBSDEs include four equations.}

A closely related work \cite{AMA20} investigated linear quadratic games with arbitrary number of players by the dynamic programming approach. For the finite number of agents, they considered the case that the population state average is known, i.e., the centralized or aggregate
sharing information structure. Here, by tackling non-standard FBSDEs we address the decentralized games where the population state average is unknown. 
In addition, there exist some works on mean-field cooperative teams with a finite number of players \cite{AM15, CA17, LFZ18}, in which the state average or 
all initial states are needed.

The main contributions of the paper are listed as follows:

\begin{itemize}
\item 		
{For the finite horizon problem, we first give necessary and sufficient conditions for the solvability of finite-population decentralized games. By decoupling non-standard FBSDEs, we design a set of decentralized Nash strategies in terms of two Riccati equations.} 
 {For the infinite horizon case, 
 a simple condition is given for the existence of the stabilizing solution to the algebraic Riccati equation arising from consensus.}
	\item   {The proposed decentralized
Nash equilibrium coincides with the result of classical mean-field games in the infinite population case. The implementation procedure of the proposed decentralized
strategies is further provided using the consensus approach.}
\item {The LQG team problem is also studied, where the social cost
is a weighted sum of individual costs. Under mild conditions, we give the decentralized social optimal control and the corresponding social cost.}
\end{itemize}

\subsection{Organization and Notation}

The organization of the paper is as follows. In Section 2, we formulate the problems of finite-population mean-field game and control. In Section 3, decentralized Nash equilibria are designed for the finite- and infinite- horizon problems, respectively. Section 4 shows the connection between the proposed decentralized control and the results of classical mean-field games. 
In Section 5, we consider the team problem. 
 Section 6 gives two numerical examples to verify results. Section 7 concludes the paper.

The following notation will be used throughout this paper. $\|\cdot\|$
denotes the Euclidean vector norm or matrix spectral norm. For a vector $z$ and a matrix $Q$, $\|z\|_Q^2= z^TQz$, and $Q>0$ ($Q\geq0$) means that $Q$ is positive definite (semi-positive definite). For two vectors $x,y$, $\langle x,y\rangle=x^Ty$.
$C([0,T],\mathbb{R}^n)$ is the space of all $\mathbb{R}^n$-valued continuous functions defined on $[0,T]$,
 and $C_{\rho/2}([0,\infty), \mathbb{R}^n)$ is a subspace of $C([0,\infty),\mathbb{R}^n)$  which is given by $\{f|\int_0^{\infty}e^{-\rho t}\|f(t)\|^2dt<\infty  \}.$

\section{Problem Description}
\label{sec2.3.1}

Consider an MAS with $N$ agents. The $i$th agent evolves by the following stochastic differential equation (SDE):
\vspace{-0.1cm}
\begin{align}\label{eq1}
dx_i(t) = \ & [Ax_i(t)+Bu_i(t)+f(t)]dt\cr
&+[Cx_i(t)+Du_i(t)+\sigma(t)] dw_i(t), 
\end{align}
where $x_i\in
\mathbb{R}^n$ and $u_i\in\mathbb{R}^r$  are the state and input of the $i$th agent, $ i=1,\cdots, N$. $f,\sigma\in C_{\rho/2}([0,\infty), \mathbb{R}^n)$ reflect
the impact on each agent by the external environment.
$\{w_i(t),i=1,\cdots, N\}$ are a sequence of independent $1$-dimensional Brownian motions on a complete
filtered probability space $(\Omega,
\mathcal F, \{\mathcal F_t\}_{ t\geq 0}, \mathbb{P})$.
The cost function of agent $i$ is given by
\vspace{-0.3cm}
\begin{equation}\label{eq2}
\begin{aligned}
J_i(u_i,u_{-i})= \mathbb{E}\int_0^{\infty}&e^{-\rho t}\Big\{\big\|x_i(t)
-\Gamma x^{(\alpha)}(t)\\&-\eta(t)\big\|^2_{Q}
+\|u_i(t)\|^2_{R}\Big\}dt,
\end{aligned}
\end{equation}
where $Q\geq 0$, $R>0$ and $\Gamma$ are constant matrices with appropriate dimension; the vector $\eta\in C_{\rho/2}([0,\infty), \mathbb{R}^n)$ and $\rho>0$ is a discount factor; $x^{(\alpha)}(t)=\sum_{j=1}^N\alpha_j^{(N)}x_j(t)$ and $u_{-i}=\{u_1,
\ldots,u_{i-1}, u_{i+1}, \ldots, u_N\}$.
 Assume that the weight allocation satisfies:

(i) $\alpha_j^{(N)}\geq 0,$ $j=1,\cdots,N$;

(ii) $\sum_{j=1}^N\alpha_j^{(N)}=1$.
\begin{remark}
The cost $J_i$ represents  the $i$th agent error of tracking an affine function of
the weighted state average. Particularly, for the case $\Gamma=I$, $\eta=0$, each agent will tend to track the weighted state average.
The latter constitutes an optimization paradigm for the consensus problem of MAS \cite{NCMH13}.
\end{remark}

%

\begin{remark}
Note that if $\alpha_j^{(N)}=\frac{1}{N}$, $j=1,\cdots,N$, then $x^{(\alpha)}(t)=x^{(N)}(t)\stackrel{\Delta}{=}\frac{1}{N}\sum_{j=1}^Nx_j(t)$.
The weighted summation $ x^{(\alpha)}$ is a generalization of the population state average $ x^{(N)}$, which was commonly considered in the classical mean-field games \cite{HCM07}, \cite{BFY13}, \cite{WZ12}. 
For this weighted cost, different agents around an individual may affect differently the individual, and some agents are allowed to be more dominant. The weighted average interaction models appear in many practical issues such as
selfish herd behavior of animals \cite{RV05}, lattice models in retailing service \cite{B93}, deep structured teams \cite{AA20} and social segregation
phenomena \cite{S71}.
See \cite{HMC10} for the NCE principle with weighted cost interactions.

\end{remark}

We start with some definitions. Denote the filtration of agent $i$ by the $\sigma$-subalgebra 
$${\mathcal F}_t^{i}=\sigma(x_{i}(0),w_i(s),0\leq s\leq t),\ i=1,\cdots,N.$$
The admissible decentralized control set is given by $$
\begin{aligned}
{\mathcal U}_{d,i} &=\Big\{u_i\ |\
u_i(t)\ \hbox{is adapted to}\  {\mathcal F}_t^{i}\Big\}.
\end{aligned}
$$

\begin{definition}
  A set of strategies $\{\check{u}_i\in {\mathcal U}_{d,i} ,i=1,\cdots,N\}$
   is said to be a decentralized Nash equilibrium if 
   the following holds:
  $$J_i(\check{u}_i,\check{u}_{-i})=\inf_{u_i\in \mathcal{U}_{d,i}}J_i(u_i,\check{u}_{-i}),\ i=1,\cdots,N.$$	
 Furthermore, $\{\check{u}_i\in {\mathcal U}_{d,i} ,i=1,\cdots,N\}$
   is  a decentralized social optimal solution if the following holds:
 $$J_{\rm soc}(\check{u})=\inf_{u_i\in \mathcal{U}_{d,i}, 1\leq i \leq N}J_{\rm soc}(u),$$	
where $J_{\rm soc}(u)=\sum_{i=1}^NJ_{i}(u)$ and $u=(u_1,\cdots,u_N)$.
\end{definition}

In this paper, we mainly study the following problems.

\textbf{(G)} Seek a set of decentralized Nash equilibrium strategies  $\{\check{u}_i\in {\mathcal U}_{d,i} ,i=1,\cdots,N\}$ with respect to cost functions
$\{J_i,i=1,\cdots,N\}$
for the system (\ref{eq1})-(\ref{eq2}).

\textbf{(S)} Seek a set of decentralized social optimal strategies  $\{\check{u}_i\in {\mathcal U}_{d,i} ,i=1,\cdots,N\}$ with respect to cost functions
$\{J_i,i=1,\cdots,N\}$
for the system (\ref{eq1})-(\ref{eq2}). 


We make the assumption on initial states of agents.

\textbf{A1)} $x_i(0)=x_{i0}, i=1,...,N$, are mutually independent and have the same mathematical expectation $\bar{x}_0$. There exists a constant $c_0$ such that $\max_{1\leq i \leq N}\mathbb{E}\|x_i(0)\|^2<c_0$.

\section{Finite-Population Mean-field LQG Games}

\subsection{The finite-horizon problem}\label{sec3a}

From now on, the time variable $t$ may be suppressed
when no confusion occurs. For the convenience of design, we first consider the following finite-horizon problem:
\begin{align*}
{\bf{(G^{\prime})}} \ \ d{x}_i=&(A{x}_i+B{u}_i
+ f)dt+
(Cx_i+Du_i+\sigma)  dw_i,
\cr
J_{i, \rm T}(u_i,u_{-i})&=\mathbb{E}\!\int_0^{T}\!\!e^{-\rho t}\big[\|x_i-\Gamma x^{(\alpha)}-\eta\big\|^2_{Q}+\|u_i\|^2_R\big]dt.
\end{align*}

We now obtain some necessary and sufficient conditions for
the existence of decentralized Nash equilibrium strategies of Problem (G$^{\prime}$)  by using variational analysis.

\begin{theorem}\label{thm1}
(i) If Problem (G$^{\prime}$) admits a set of decentralized Nash equilibrium strategies $\{\check{u}_i\in \mathcal{U}_{d,i}, i=1,\cdots,N\}$, then the following FBSDE system admits a set of adapted solutions $(\check{x}_i,\check{\lambda}_i, \check{\beta}_{i}^{j},i,j=1,\cdots,N)$:
\begin{equation}\label{eq6a}
\left\{  \begin{aligned}
d\check{x}_i=&(A\check{x}_i+B\check{u}_i+f)dt
  +(C\check{x}_i+D\check{u}_i+\sigma) dw_i,\\
d\check{\lambda}_i=&-\Big[(A-\rho I)^T\check{\lambda}_i+C^T\check{\beta}_i^i\Big]dt+ \sum_{j=1}^N
\check{\beta}_{i}^j dw_j\cr&+\Big[\Big(I-\alpha_i^{(N)}\Gamma\Big)^T Q(\check{x}_i-\Gamma \check{x}^{(\alpha)}-\eta)\Big]dt,\cr
\check{x}_i(0)&=x_{i0},\quad  \check{\lambda}_i(T)=0, \ i=1,\cdots,N
\end{aligned}\right.
\end{equation}
with
\begin{equation}\label{eq7a}
	  \check{u}_i=-R^{-1}\big(B^T\mathbb{E}[\check{\lambda}_i|\mathcal{F}_t^i]+D^T\mathbb{E}[\check{\beta}^i_i|\mathcal{F}_t^i]\big).
  \end{equation}
      (ii) If the equation system (\ref{eq6a}) admits a set of solutions $(\check{x}_i,\check{\lambda}_i, \check{\beta}_{i}^{j},i,j=1,\cdots,N)$, then Problem (G$^{\prime}$) has a set of decentralized Nash equilibrium strategies $\check{u}_i, i=1,\cdots,N$, which is given by (\ref{eq7a}). 
\end{theorem}
\emph{Proof.} 
See Appendix A. \hfill{$\Box$}

\begin{remark}
The definition of adapted solutions arises in solving backward SDEs. The adapted solution to the backward SDE in (\ref{eq6a}) is a
sequence of adapted stochastic processes  $(\check{\lambda}_i, \check{\beta}_{i}^{j},i,j=1,\cdots,N)$, and it is the terms $\check{\beta}_i^j, j=1,\cdots,N$ that correct the possible ``nonadaptiveness" caused by the backward nature of 
(\ref{eq6a}).  See \cite{MY99} for more details of adapted solutions.
\end{remark}

\subsubsection{Homogeneous weights}
We first consider the case $\alpha_i^{(N)}=\frac{1}{N}$.
{
\begin{lemma}\label{lem1a}
For any $j\not=i$, the following holds:
\begin{equation}\label{eq3}
\begin{split}
		\mathbb{E}\big[\check{x}_j(r)|\mathcal{F}_t^i\big]&=\mathbb{E}\big[\check{x}_j(r)\big]=
		\mathbb{E}\big[\check{x}_i(r)\big].
	\end{split}
\end{equation}
\end{lemma}
\emph{Proof.}
Note that $\check{x}_i(t)$ is adapted to $\mathcal{F}_t^i$, $i=1,\cdots,N$. Since 
$w_i, i=1,...,N$, are mutually independent, by A1) then $\check{x}_i, i=1,...,N$, are independent of each other. Since all agents have the same parameters, we obtain (\ref{eq3}). \hfill{$\Box$}
}

By Lemma \ref{lem1a},
 we have
\begin{equation}\nonumber
	\begin{split} \mathbb{E}[\check{x}^{(N)}|\mathcal{F}_t^i]&=\frac{1}{N}\check{x}_i+\frac{1}{N}\sum_{j\not=i}\mathbb{E}[\check{x}_j]=\frac{1}{N}\check{x}_i+\frac{N-1}{N}\mathbb{E}[\check{x}_i].
	\end{split}
\end{equation}
It follows from (\ref{eq6a})-(\ref{eq7a}) that
\begin{equation}\label{eq11}
\left\{
\begin{aligned}
d\mathbb{E}[\check{x}_i]=&(A\mathbb{E}[\check{x}_i]-BR^{-1} (B^T\mathbb{E}[\check{\lambda}_i]+D^T[\check{\beta}_i^i]+f)dt, \cr
    d\mathbb{E}[\check{\lambda}_i|\mathcal{F}_t^i]&= - \Big[(A-\rho I)^T\mathbb{E}[\check{\lambda}_i|\mathcal{F}_t^i]+(I-\frac{1}{N}\Gamma)^T Q\cr&\times\big((I-\frac{1}{N}\Gamma)\check{x}_i-\frac{N-1}{N}\Gamma\mathbb{E}[\check{x}_i]-\eta\big) \Big]dt\cr
    &-C^T\mathbb{E}[\check{\beta}^i_i|\mathcal{F}_t^i]dt+\mathbb{E}[\check{\beta}^i_i|\mathcal{F}_t^i]dw_i,
    \\ \mathbb{E}[\check{x}_i(0)]&=\bar{x}_{0},\quad \mathbb{E}[\check{\lambda}_i(T)|\mathcal{F}_T^i]=0.
\end{aligned}\right.
\end{equation}

\begin{proposition}\label{prop3.0}
The FBSDE (\ref{eq6a}) admits a set of adapted solutions 
if and only if FBSDE
   (\ref{eq11}) admits a set of adapted solutions. 
\end{proposition}
\emph{Proof.} If
   (\ref{eq11}) admits an adapted solution, then $(\check{u}_i,i=1,\cdots,N)$ is given and (\ref{eq6a}) is decoupled. Then (\ref{eq6a}) admits a solution.
{
   The remainder of the proof is straightforward.} \hfill$\Box$


We now study under what conditions the decentralized Nash equilibrium (\ref{eq7a}) has a feedback representation.
Let $K_N,\Pi_N\in \mathbb{R}^{n\times n}$ and $s_N\in \mathbb{R}^n$ satisfy
\begin{align}\label{eq9a}
\rho K_N=&\dot{K}_N+{A}^TK_N+K_NA-(B^TK_N+D^TK_NC)^T\cr&\times\Upsilon_N^{-1}(B^TK_N+D^TK_NC)+C^TK_NC\cr&+\Big(I-\frac{1}{N}\Gamma\Big)^T Q\Big(I-\frac{1}{N}\Gamma\Big), \ K_N(T)=0,
\\
 \label{eq9b}
\rho \Pi_N=& \dot{\Pi}_N+{A}^T\Pi_N+\Pi_N{A}-\Pi_N B \Upsilon_N^{-1}B^T\Pi_N\cr&-\Pi_N B \Upsilon_N^{-1}(B^TK_N+D^TK_NC)
\cr&-(B^TK_N+D^TK_NC)^T\Upsilon_N^{-1}B^T\Pi_N\cr&-\frac{N-1}{N}\Big(I-\frac{1}{N}\Gamma\Big)^T Q\Gamma,\
 \Pi_N(T)=0,\\
\label{eq10}
\rho s_N\! =&\dot{s}_N\!+\!\big[{A}\!-\!B\Upsilon_N^{-1}\big( B^T(K_N\!+\!\Pi_N)\!+\!D^TK_NC\big)\big]^Ts_N\cr&+(K_N+\Pi_N)f-(I-\frac{1}{N}\Gamma)^TQ\eta+\big[C-D\Upsilon_N^{-1}\cr&\cdot\big( B^T(K_N\!+\!\Pi_N)\!+\!D^TK_NC\big)\big]^T\!K_N\sigma, s_N(T)=0,\cr
\end{align}
with $\Upsilon_N\stackrel{\Delta}{=}R+D^TK_ND.$
\begin{lemma}\label{lem1}
  If (\ref{eq9a})-(\ref{eq10}) have a solution, respectively, then (\ref{eq6a}) admits a set of adapted solutions.
\end{lemma}
\emph{Proof.}  Let
\begin{equation}\label{eq12b}
 \bar{s}_N(t)= \mathbb{E}[\check{\lambda}_i(t)|\mathcal{F}_t^i]-K_N(t)\check{x}_i(t)-\Pi_N(t)\mathbb{E}[\check{x}_i(t)],
\end{equation}
   where $K_N(t)$ and $\Pi_N(t)$ satisfy (\ref{eq9a}) and (\ref{eq9b}), respectively. Here, $\bar{s}_N(t)\in \mathbb{R}^n$ may depend on $i$. However, we will show later that $\bar{s}_N(t)$ is actually independent of $i$. It follows by (\ref{eq12b}) that $\mathbb{E}[\check{\lambda}_i(t)]=(K_N(t)+\Pi_N(t))\mathbb{E}[\check{x}_i(t)]+\bar{s}_N(t)$.
 By It\^{o}'s formula, we obtain
\begin{align}\label{eq13b}
d\mathbb{E}[\check{\lambda}_i|\mathcal{F}_t^i]\!=&\dot{K}_N\check{x}_idt+K_N\big[(A\check{x}_i+B\check{u}_i+f)dt\cr&
\!\!+\!(C\check{x}_i+D\check{u}_i+\sigma) dw_i\big]+\dot{\Pi}_N\mathbb{E}[\check{x}_i]dt\cr&
\!\!+\!\Pi_N\big[A\mathbb{E}[\check{x}_i]+B \mathbb{E}[\check{u}_i]\!+\!f\big]dt\!+\!\dot{\bar{s}}_Ndt.
\end{align}
Comparing this with (\ref{eq11}) and equating the $dw_i$ terms, it follows that $\mathbb{E}[\check{\beta}^i_i|\mathcal{F}_t^i]=K_N
(C\check{x}_i+D\check{u}_i+\sigma)$.
By (\ref{eq7a}) and (\ref{eq12b}), we have
\begin{equation}\nonumber
	\begin{split}
		&R\check{u}_i+B^T(K_N\check{x}_i+\Pi_N\mathbb{E}[\check{x}_i]+\bar{s}_N)\\&+D^TK_N
		(C\check{x}_i+D\check{u}_i+\sigma)=0,
	\end{split}
\end{equation}
which leads to
\begin{equation}\label{eq15e}
	\begin{split}
	  \check{u}_i=&-(R+D^TK_ND)^{-1}\big[(B^TK_N+D^TK_NC)\check{x}_i\cr&+B^T\Pi_N\mathbb{E}[\check{x}_i]+B^T\bar{s}_N+D^TK_N\sigma\big].
	\end{split}
  \end{equation}
 Then
\begin{equation}\label{eq8}
	\begin{split}
  \mathbb{E}[\check{u}_i]=&-\Upsilon_N^{-1}\big[(B^TK_N+B^T\Pi_N+D^TK_NC)\mathbb{E}[\check{x}_i]\cr&+B^T\bar{s}_N+D^TK_N\sigma\big].
	\end{split}
  \end{equation}
Applying (\ref{eq15e})-(\ref{eq8}) to (\ref{eq11}) and (\ref{eq13b}), and comparing them, we obtain
that $\bar{s}_N$ satisfies (\ref{eq10}). Here, the $\check{x}_i$ terms equaling 0 is obtained from (\ref{eq9a}), and  the $\mathbb{E}[\check{x}_i]$ terms equalling 0 from (\ref{eq9b}).
From (\ref{eq12b}), FBSDE (\ref{eq11}) is solvable. By Proposition \ref{prop3.0}, (\ref{eq6a}) admits a set of adapted solutions.
$\hfill \Box$

By Lemma \ref{lem1} and (\ref{eq15e}), we obtain the following decentralized strategies for $N$ agents,
\begin{equation}\label{eq15a}
	\begin{split}
	 \check{u}_i=&-\Upsilon_N^{-1}\big[(B^TK_N+D^TK_NC)\check{x}_i+B^T\Pi_N\mathbb{E}[\check{x}_i]\cr&+B^Ts_N+D^TK_N\sigma\big], \quad i=1,\cdots, N,
	\end{split}
\end{equation}
where $K_N, \Pi_N, s_N$ are given by (\ref{eq9a})-
(\ref{eq10}),  and $\mathbb{E}[\check{x}_i]$ obeys
\begin{align}\label{eq19}
  d\mathbb{E}[\check{x}_i]\!=&\big[A-B\Upsilon^{-1}_N\big(B^TK_N+B^T\Pi_N+D^TK_NC\big)\big]\cr
   &\times\mathbb{E}[\check{x}_i]dt-B\Upsilon_N^{-1}(B^T s_N+D^TK_N\sigma)dt+fdt,\cr
   \mathbb{E}[\check{x}_i(0)]&=\bar{x}_{0}.
  \end{align}

\begin{remark}
 In previous works (e.g., \cite{HCM07}, \cite{LZ08}),
the mean-field term $x^{(N)}$ in cost functions is first substituted by a deterministic function $\bar{x}$. By handling the fixed-point equation, 
 $\bar{x}$ is obtained, and the decentralized control is constructed. 
 Here, we first obtain the solvability condition of decentralized control by virtue of variational analysis, and then design decentralized control laws
  by tackling FBSDEs and conditional Hamiltonians. Note that in this case $s_N$ and $\mathbb{E}[\check{x}_i]$ are decoupled and no fixed-point equation is needed.
\end{remark}

Next, we study whether that above proposed strategy gives a decentralized Nash equilibrium.
For the analysis, we introduce the following assumption.

\textbf{A2)} Equation (\ref{eq9b}) admits a solution in $C([0,T], \mathbb{R}^{n\times n})$.

\begin{theorem} \label{thm3}
	Let A1)-A2) hold.
Then for Problem (G$^{\prime}$), the set of control laws
	$\{\check{u}_1,\cdots,\check{u}_N\}$ given by (\ref{eq15a}) is a decentralized Nash equilibrium, 
 i.e.,
  $$J_{i,\rm T}(\check{u}_i,\check{u}_{-i})=\inf_{u_i\in \mathcal{U}_{d,i}}J_{i,\rm T}(u_i,\check{u}_{-i}).$$	
\end{theorem}

\emph{Proof.} Note that (\ref{eq9a}) admits a solution in $C([0,T], \mathbb{R}^{n\times n})$ since $\Big(I-\frac{1}{N}\Gamma\Big)^T Q\Big(I-\frac{1}{N}\Gamma\Big)\geq 0$. Under A2), (\ref{eq10}) has a solution in $C([0,T], \mathbb{R}^{n})$. By Lemma \ref{lem1}, (\ref{eq6a}) admits a set of adapted solutions,
which together with Theorem \ref{thm1} completes the proof of the theorem. $\hfill \Box$

\begin{remark}
	In this paper, we consider the case $Q\geq 0$ and $R>0$. Indeed, even if $Q$ and $R$ are indefinite, under some convex conditions, we may design the decentralized strategies and obtain the optimality result.
	(See e.g., \cite{SLY16})
\end{remark}

By \cite[P. 48]{MY99}, A2) holds if and only if $$\det\Big\{[0, \ I]e^{\mathcal{A}t}\left[
\begin{array}{c}
  0\\
 I
\end{array}
\right]\Big\}>0,\qquad \forall t\in [0,T],$$
where $	\bar{A}_N\stackrel{\Delta}{=}A-B \Upsilon_N^{-1}(B^TK_N+D^TK_NC)$ and $$\mathcal{A}\stackrel{\Delta}{=}
		 {\Bigg[\begin{array}{cc}
			\bar{A}_N-\frac{\rho}{2}I&-B \Upsilon_N^{-1}B^T\\
				\frac{N-1}{N}\Big(I-\frac{1}{N}\Gamma\Big)^T Q\Gamma&-(\bar{A}_N-\frac{\rho}{2}I)^T
			\end{array}
			\Bigg].}
$$
Particularly, if $\Gamma=I$, we have the following result.

\begin{proposition}\label{prop3.1}
  If $\Gamma=I$, then A2) holds necessarily.
\end{proposition}
\emph{Proof.}  See Appendix A. \hfill{$\Box$}

\subsubsection{Heterogeneous weights}

We now consider the situation that some agent plays a dominant role in the system. Specifically,
the weights in $x^{(\alpha)}$ are taken as
$\alpha_1=\alpha$, and $\alpha_i=\frac{1-\alpha}{N-1}$, $ i=2,\cdots, N$, $0< \alpha < 1$.
For simplicity, consider the case where $f=\sigma=\eta=0$.

Note that
$\mathbb{E}\big[\check{x}_2\big]=\cdots=\mathbb{E}\big[\check{x}_N\big]\not=\mathbb{E}\big[\check{x}_1\big]$.
By Lemma \ref{lem1a}, we obtain that for $i=1$,
$\mathbb{E}[\check{x}^{(\alpha)}|\mathcal{F}_t^1]=\alpha \check{x}_1+(1-\alpha)\mathbb{E}[\check{x}_j]\ (j\not=1),$
and
$i\not=1$, $$\mathbb{E}[\check{x}^{(\alpha)}|\mathcal{F}_t^i]=\frac{1-\alpha}{N-1}\check{x}_i+\alpha \mathbb{E}[\check{x}_1]+\frac{N-2}{N-1}(1-\alpha)\mathbb{E}[\check{x}_i].$$
It follows from (\ref{eq6a}) that
\begin{equation}\label{eq11c}
\left\{
\begin{split}
&d\mathbb{E}[\check{x}_i]=(A\mathbb{E}[\check{x}_i]+B \mathbb{E}[\check{u}_i])dt, \mathbb{E}[\check{x}_i(0)]=\bar{x}_{0},\\
&\qquad\qquad\qquad i=1,\cdots, N,\cr
    &d\mathbb{E}[\check{\lambda}_1|\mathcal{F}_t^1]= - \big[(A-\rho I)^T\mathbb{E}[\check{\lambda}_1|\mathcal{F}_t^1]
    +C^T\mathbb{E}[\check{\beta}^1_1|\mathcal{F}_t^1]dt\\&+(I-\alpha\Gamma)^T Q\big((I-\alpha\Gamma)\check{x}_1-(1-\alpha)\Gamma\mathbb{E}[\check{x}_j]\big)\big]dt\cr&+\mathbb{E}[\check{\beta}^1_1|\mathcal{F}_t^1]dw_1, \quad \mathbb{E}[\check{\lambda}_1(T)|\mathcal{F}_T^1]=0, \cr
    &d\mathbb{E}[\check{\lambda}_j|\mathcal{F}_t^j]= - \Big[(A-\rho I)^T\mathbb{E}[\check{\lambda}_j|\mathcal{F}_t^j]+C^T\mathbb{E}[\check{\beta}^j_j|\mathcal{F}_t^j] \\&+\big(I-\frac{1-\alpha}{N-1}\Gamma\big)^T\times Q\Big(\big(I-\frac{1-\alpha}{N-1}\Gamma\big)\check{x}_j-\alpha \Gamma\mathbb{E}[x_1]\cr&-\frac{N-2}{N-1}(1-\alpha)\Gamma\mathbb{E}[\check{x}_j]\Big) \Big]dt+\mathbb{E}[\check{\beta}^j_j|\mathcal{F}_t^j]dw_j, \ \cr&\mathbb{E}[\check{\lambda}_j(T)|\mathcal{F}_T^j]=0,\quad j=2,\cdots,N.
\end{split}\right.
\end{equation}

We now
 decouple 
(\ref{eq11c}) by using the idea of Lemma \ref{lem1} (see also \cite{MY99, QZW19}).
Let $$\mathbb{E}[\check{\lambda}_1|\mathcal{F}_t^1]=\bar{K}^1_N\check{x}_1+\bar{\Pi}^{1,1}_N \mathbb{E}[\check{x}_1]+\bar{\Pi}^{1,j}_N \mathbb{E}[\check{x}_j], \ j\not= 1,$$
where $\bar{K}^1_N,\bar{\Pi}^{1,1}_N , \bar{\Pi}^{1,j}_N \in \mathbb{R}^{n\times n} $. Then by It\^{o}'s formula, 
\begin{align}\label{eq29}
d\mathbb{E}[\check{\lambda}_1|\mathcal{F}_t^1]=&\dot{\bar{K}}^1_N\check{x}_1dt+\bar{K}^1_N\big[(A\check{x}_1+B\check{u}_1)dt\cr&
+(C\check{x}_1+D\check{u}_1) dw_1\big]+\big[\dot{\bar{\Pi}}_N^{1,1}\mathbb{E}[\check{x}_1]\cr&
+\bar{\Pi}_N^{1,1}(A\mathbb{E}[\check{x}_1]+B\mathbb{E}[\check{u}_1])+\dot{\bar{\Pi}}_N^{1,j}\mathbb{E}[\check{x}_j]\big]dt
  \cr&+\bar{\Pi}_N^{1,j}(A\mathbb{E}[\check{x}_j]+B\mathbb{E}[\check{u}_j])dt.\cr
\end{align}
Comparing this with (\ref{eq11c}), it follows that $\mathbb{E}[\check{\beta}^1_1|\mathcal{F}_t^i]=\bar{K}_N^1
(C\check{x}_1+D\check{u}_1)$. By (\ref{eq7a}), we have
\begin{equation}\nonumber
	\begin{split}
	 R\check{u}_i+B^T(\bar{K}^1_N\check{x}_1+\bar{\Pi}^{1,1}_N \mathbb{E}[\check{x}_1]+\bar{\Pi}^{1,j}_N \mathbb{E}[\check{x}_j])\\+D^T\bar{K}^1_N
	(C\check{x}_1+D\check{u}_1)=0.	
	\end{split}
\end{equation}
This leads to
\begin{align}\label{eq30}
  \check{u}_1=&-(\bar{\Upsilon}_N^1)^{-1}\big[(B^T\bar{K}^1_N+D^T\bar{K}^1_NC)\check{x}_1\cr&+B^T\bar{\Pi}^{1,1}_N \mathbb{E}[\check{x}_1]+B^T\bar{\Pi}^{1,j}_N \mathbb{E}[\check{x}_j]\big], \quad j\not= 1,
  \end{align}
where $\bar{\Upsilon}_N^1\stackrel{\Delta}{=}R+D^T\bar{K}^1_ND.$

Let $\mathbb{E}[\check{\lambda}_j|\mathcal{F}_t^j]=\bar{K}^j_N\check{x}_j+\bar{\Pi}^{j,1}_N \mathbb{E}[\check{x}_1]+\bar{\Pi}^{j,j}_N \mathbb{E}[\check{x}_j]$, $2\leq j\leq N$,  where $\bar{K}^j_N,\bar{\Pi}^{j,1}_N , \bar{\Pi}^{j,j}_N \in \mathbb{R}^{n\times n} $. Then by It\^{o}'s formula, 
\begin{align}\label{eq31a}
&d\mathbb{E}[\check{\lambda}_j|\mathcal{F}_t^j]\cr=&\dot{\bar{K}}^j_N\check{x}_jdt+\bar{K}^j_N\big[(A\check{x}_j+B\check{u}_j)dt
+(C\check{x}_j+D\check{u}_j) dw_j\big]\cr&+\dot{\bar{\Pi}}_N^{j,1}\mathbb{E}[\check{x}_1]dt
+\bar{\Pi}_N^{j,1}(A\mathbb{E}[\check{x}_1]+B\mathbb{E}[\check{u}_1])dt\cr
& +\dot{\bar{\Pi}}_N^{j,j}\mathbb{E}[\check{x}_j]dt
+\bar{\Pi}_N^{j,j}(A\mathbb{E}[\check{x}_j]+B\mathbb{E}[\check{u}_j])dt.
\end{align}
Comparing this with (\ref{eq11c}), it follows that $\mathbb{E}[\check{\beta}^j_j|\mathcal{F}_t^j]=\bar{K}_N^j
(C\check{x}_j+D\check{u}_j)$. This together with (\ref{eq7}) gives
\begin{equation}
	\begin{split}
&		R\check{u}_j+B^T(\bar{K}^j_N\check{x}_j+\bar{\Pi}^{j,1}_N \mathbb{E}[\check{x}_1]+\bar{\Pi}^{j,j}_N \mathbb{E}[\check{x}_j])\\
&+D^T\bar{K}^j_N
		(C\check{x}_j+D\check{u}_j)=0.
	\end{split}
\end{equation}
This leads to
\begin{align}\label{eq32a}
  \check{u}_j=&-(\bar{\Upsilon}_N^j)^{-1}\big[(B^T\bar{K}^j_N+D^T\bar{K}^j_NC)\check{x}_j\cr&+B^T\bar{\Pi}^{j,1}_N \mathbb{E}[\check{x}_1]+B^T\bar{\Pi}^{j,j}_N \mathbb{E}[\check{x}_j]\big], 
  \end{align}
   where $\bar{\Upsilon}_N^j\stackrel{\Delta}{=}R+D^T\bar{K}^j_ND, 2\leq j\leq N.$

Applying (\ref{eq30}) and (\ref{eq32a}) to (\ref{eq11c})-(\ref{eq29}) and comparing them, it follows that
\begin{align}\label{eq37a}
\rho\bar{K}_N^i\!=&\dot{\bar{K}}_N^1+{A}^T\bar{K}^1_N+\bar{K}^1_NA-(B^T\bar{K}^1_N+D^T\bar{K}^1_NC)^T\cr&\times(\bar{\Upsilon}_N^1)^{-1}(B^T\bar{K}^1_N+D^T\bar{K}^1_NC)+C^T\bar{K}^1_NC\cr&+\big(I-\alpha\Gamma\big)^T Q\big(I-\alpha\Gamma\big), \quad\bar{K}^1_N(T)=0,\\
\label{eq37b}\rho\bar{\Pi}^{1,1}_N\!=&\dot{\bar{\Pi}}^{1,1}_N\!+{A}^T\bar{\Pi}^{1,1}_N\!+\bar{\Pi}^{1,1}_N{A}\!-\bar{\Pi}_N^{1,1} B (\bar{\Upsilon}_{ N}^j)^{-1}B^T\bar{\Pi}_N^{1,1}\cr&-(B^T\bar{K}^j_N+D^T\bar{K}^j_NC)^T(\bar{\Upsilon}_{ N}^1)^{-1}B^T\bar{\Pi}_N^{1,1}
\cr&-\bar{\Pi}_N^{j,1} B (\bar{\Upsilon}_{ N}^1)^{-1}(B^T\bar{K}^1_N+D^T\bar{K}^1_NC\!+\!B^T\bar{\Pi}_N^{1,1})\cr&-\bar{\Pi}^{1,j}_N B (\bar{\Upsilon}_{ N}^j)^{-1}B^T\bar{\Pi}^{j,1}_N,
\ \quad\bar{\Pi}_N^{1,1}(T)=0,\\
\label{eq37c}
\rho\bar{\Pi}^{1,j}_N\!= &\dot{\bar{\Pi}}^{1,j}_N+{A}^T\bar{\Pi}^{1,j}_N+\bar{\Pi}^{1,j}_N{A}-(B^T\bar{K}^1_N+D^T\bar{K}^1_NC\cr&
+B^T\bar{\Pi}^{1,1}_N)^T (\bar{\Upsilon}_{ N}^1)^{-1}B^T\bar{\Pi}_N^{1,j}-\bar{\Pi}_N^{1,j}B(\bar{\Upsilon}_{ N}^j)^{-1}\cr
&\times(B^T\bar{K}^j_N+D^T\bar{K}^j_NC+B^T\bar{\Pi}^{j,j}_N)\cr&
-(1-\alpha)(I-\alpha\Gamma)^T Q\Gamma,
\quad \bar{\Pi}_N^{1,j}(T)=0.
\end{align}
In the above,  (\ref{eq37a}) is obtained by equating the $\check{x}_1$ terms; (\ref{eq37b}) by equating the $\mathbb{E}[\check{x}_1]$ terms and (\ref{eq37c}) by equating the $\mathbb{E}[\check{x}_j]$ terms.
After strategies (\ref{eq30}) and (\ref{eq32a}) are applied, comparing (\ref{eq11c}) and (\ref{eq31a}), it follows that
\begin{align}\label{eq41a}
\rho\bar{K}^{j}_N\!=&\dot{\bar{K}}_N^j+{A}^T\bar{K}^j_N+\bar{K}^j_NA-(B^T\bar{K}^j_N+D^T\bar{K}^j_NC)^T
\cr\times&(\bar{\Upsilon}_N^j)^{-1}(B^T\bar{K}^j_N+D^T\bar{K}^j_NC)+C^T\bar{K}^j_NC\cr+&\big(I-\frac{1\!-\!\alpha}{N-1}\Gamma\big)^T Q\big(I\!-\!\frac{1-\alpha}{N-1}\Gamma\big)=0,  \bar{K}^j_N(T)=0,
\\
\label{eq41b}
\rho\bar{\Pi}^{j,1}_N\!= &\dot{\bar{\Pi}}^{j,1}_N+{A}^T\bar{\Pi}^{j,1}_N+\bar{\Pi}^{j,1}_N{A}-\bar{\Pi}^{j,j}_N B (\bar{\Upsilon}_{ N}^1)^{-1}B^T\bar{\Pi}^{j,1}_N\cr&-(B^T\bar{K}^1_N+D^T\bar{K}^1_NC)^T(\bar{\Upsilon}_{ N}^j)^{-1}B^T\bar{\Pi}^{j,1}_N\cr&-\bar{\Pi}_N^{j,j}B(\bar{\Upsilon}_{ N}^1)^{-1}(B^T\bar{K}^1_N+D^T\bar{K}^1_NC)\cr&-\alpha\big(I-\frac{1-\alpha}{N-1}\Gamma\big)^T Q \Gamma=0, \ \bar{\Pi}_N^{j,1}(T)=0,
\end{align}
\begin{align}
\label{eq41c}
\rho\bar{\Pi}^{j,j}_N\!=&\dot{\bar{\Pi}}^{j,j}_N+{A}^T\bar{\Pi}^{j,j}_N+\bar{\Pi}^{j,j}_N{A}\cr
-&(B^T\bar{K}^j_N+D^T\bar{K}^j_NC)^T(\bar{\Upsilon}_{ N}^j)^{-1}B^T\bar{\Pi}^{j,j}_N\cr
-&\bar{\Pi}_N^{j,j}B(\bar{\Upsilon}_{ N}^j)^{-1}(B^T\bar{K}^j_N+D^T\bar{K}^j_NC)
\cr-&\bar{\Pi}^{j,j}_N B (\bar{\Upsilon}_{ N}^j)^{-1}B^T\bar{\Pi}^{j,j}_N-\bar{\Pi}^{j,1}_N B (\bar{\Upsilon}_{ N}^1)^{-1}B^T\bar{\Pi}^{1,j}_N\cr
-&\frac{N-2}{N-1}(1-\alpha)\big(I-\frac{1-\alpha}{N-1}\Gamma\big)^T Q,
\ \bar{\Pi}_N^{j,j}(T)=0.\cr
\end{align}

\begin{theorem} Assume that A1) holds. For Problem (G$^{\prime}$), if (\ref{eq37b})-(\ref{eq37c}) and (\ref{eq41b})-(\ref{eq41c}) admit solutions, respectively, then the set of strategies given by (\ref{eq30}) and (\ref{eq32a}) is a decentralized Nash equilibrium.
\end{theorem}
\emph{Proof.}
 Note that $\big(I-\alpha\Gamma\big)^T Q\big(I-\alpha\Gamma\big)\geq 0$ and $\big(I-\frac{1-\alpha}{N-1}\Gamma\big)^T Q\big(I-\frac{1-\alpha}{N-1}\Gamma\big)\geq 0$. Then (\ref{eq37a}) and (\ref{eq41a}) admit solutions $\bar{K}_N^1\geq 0$ and $\bar{K}_{N}^j\geq0$, respectively. If (\ref{eq37b})-(\ref{eq37c}) and (\ref{eq41b})-(\ref{eq41c}) admit solutions, respectively, then by the derivation above, we obtain that (\ref{eq6a}) admits a set of adapted solutions. By Theorem \ref{thm1},
the set of strategies in (\ref{eq30}) and (\ref{eq32a}) is a decentralized Nash equilibrium. $\hfill \Box$


\subsection{The infinite-horizon problem}\label{sec4}
In this section, we consider the infinite-horizon problem with homogeneous weights $(\alpha_i^{(N)}=\frac{1}{N})$.
Based on the analysis in Section \ref{sec3a}, we may design the following decentralized control for Problem (G):
\begin{equation}\label{eq14c}
\begin{aligned}
\check{u}_i(t)=&-\Upsilon_N^{-1}\big[(B^TK_N+D^TK_NC)\check{x}_i(t)\cr&+B^T(P_N-K_N)\mathbb{E}[\check{x}_i(t)]+B^Ts_N(t)\\
&+D^TK_N\sigma(t)\big],\quad  t\geq 0, \ i=1,\cdots, N,
\end{aligned}
\end{equation}
where $K_N$ and $P_N$ satisfy 
\begin{align}\label{eq15c}
\rho K_N=&A^TK_N+K_NA-(B^TK_N+D^TK_NC)^T\Upsilon_N^{-1}\cr&\times(B^TK_N+D^TK_NC)+C^TK_NC\cr&+\big(I-\frac{1}{N}\Gamma\big)^T Q\big(I-\frac{1}{N}\Gamma\big),\\
\label{eq16c}
\rho P_N=&{A}^TP_N+P_N{A}-(B^TP_N+D^TK_NC)^T \Upsilon_N^{-1}\cr&\times(B^TP_N+D^TK_NC)+C^TK_NC\cr&+\big(I-\frac{1}{N}\Gamma\big)^T Q(I-\Gamma),
\end{align}
and
${s}_N\in C_{\rho/2}([0,\infty),\mathbb{R}^n)$ is determined by
\begin{align}
\rho s_N=&\dot{s}_N+[A-B\Upsilon_N^{-1}(B^TP_N+D^TK_NC)]^T{s}_N\cr&+\big[C-D\Upsilon_N^{-1}\big( B^TP_N+D^TK_NC\big)\big]^TK_N\sigma\cr&
+P_N f-(I-\frac{1}{N}\Gamma)^TQ\eta . \label{eq17c}
\end{align}
Here the existence conditions of $K_N, \Pi_N$ and $s_N$ need to be further investigated.

%
%

We now introduce the following assumptions, where the definitions of stabilizability and detectability can be referred in e.g., \cite{ZZC08}, \cite{QZW19}.

\textbf{A3)} The system $[A-\frac{\rho}{2}I, B; C,D]$ is stabilizable, and $[A-\frac{\rho}{2}I, C,\sqrt{Q}]$ is exactly detectable. 

\textbf{A4)} Assume that (\ref{eq16c}) admits a $\rho$-stabilizing solution for all $N>1$. 

Denote $$
\begin{aligned}
  \hat{A}_N&\stackrel{\Delta}{=}A-\frac{\rho}{2}I -B\Upsilon_N^{-1}D^TK_NC,\\
\hat{C}_N&\stackrel{\Delta}{=}C^TK_NC-C^TK_ND\Upsilon_N^{-1}D^TK_NC.
\end{aligned}$$
By \cite[Theorem 18]{HZ20},  A4) holds if and only if
$${\mathcal{M}_{\Gamma}\stackrel{\Delta}{=}
		\left[\begin{array}{cc}
			\hat{A}_N-\frac{\rho}{2}I&B \Upsilon_N^{-1}B^T\\
			\hat{C}_N+\big(I-\frac{1}{N}\Gamma\big)^T Q(I-\Gamma)  &-\hat{A}_N^T+\frac{\rho}{2}I
		\end{array}
		\right]}$$
 is $(n,n)$ c-splitting, i.e., both the open left plane and  the open right plane contain $n$ eigenvalues, respectively. 
Particularly, if $\Gamma=I$, we have the following result.  Denote
$$ \mathcal{M}_{I}\stackrel{\Delta}{=}
			\left[\begin{array}{cc}
			\hat{A}_N-\frac{\rho}{2}I&B \Upsilon_N^{-1}B^T\\
			\hat{C}_N &-\hat{A}_N^T+\frac{\rho}{2}I
		\end{array}
		\right].$$
\begin{proposition}\label{lem2c}
	For the case $\Gamma=I$, let A3) hold. Then Assumption A4) holds if and only if $M_{I}$ has no eigenvalues on the imaginary axis.
Furthermore, if $C=0$, the necessary and sufficient condition ensuring A4) is that $A-\frac{\rho}{2}I$ has no eigenvalues on the imaginary axis.
\end{proposition}

\emph{Proof.} See Appendix B.
$\hfill \Box$

\begin{example}
 Consider a two-dimensional system with $A=\left[\begin{array}{cc}
a& b \\
c& d
\end{array}\right]$, $C=0$ and $\Gamma=I$. Denote $\bar{a}=a-\frac{\rho}{2}$, and $\bar{d}=d-\frac{\rho}{2}$. We have
$  \big|\lambda I-(A-\frac{\rho}{2}I)\big|=\lambda^2-(\bar{a}+\bar{d})\lambda+\bar{a}\bar{d}-bc.$
 Then $A-\frac{\rho}{2}I$ has eigenvalues on the imaginary axis if and only if $\bar{a}+\bar{d}=0$ and $\bar{a}\bar{d}-bc\geq0$. Thus, if $\bar{a}+\bar{d}\not=0$ or $\bar{a}\bar{d}-bc<0$, then $A-\frac{\rho}{2}I$ has no eigenvalues on the imaginary axis, which with $\Gamma=I$ and $C=0$ implies that A4) holds. Particularly, if $A=0$, then A4) holds.
\end{example}


The next theorem characterizes the performance of the decentralized strategies.
\begin{theorem}\label{thm11}
	Assume A1), A3)-A4) hold and $N$ is sufficiently large such that $I-\frac{1}{N}\Gamma$ is nonsingular. For Problem (G), the set of strategies
	$\{\check{u}_1,\cdots,\check{u}_N\}$ given by (\ref{eq14c}) is a decentralized Nash equilibrium, i.e., for any $i=1,\cdots,N$,
	$J_i(\check{u}_i,\check{u}_{-i})=\inf_{u_i\in \mathcal{U}_{d,i}}J_i(u_i,\check{u}_{-i}).$
	
\end{theorem}

\emph{Proof.} See Appendix B. $\hfill \Box$

\subsubsection{The model of noisy multiple integrators}
For the case $A=C=0$, $\Gamma=B=D=I$, and $f=\eta=0$, the system (\ref{eq1})-(\ref{eq2}) reduces to the model of noisy multiple integrators. Specifically, agent $i$ evolves by
\begin{equation}\label{eq61}
  dx_i=u_idt+u_idw_i, \ i=1,\cdots,N,
\end{equation}
and the cost function is given by
\begin{align}\label{eq62}
{J}_i
=&
\mathbb{E}\!\int_0^{\infty}\!e^{-\rho t}\big\{\big\|x_i
- {x}^{(N)}\big\|_Q^2
+\|u_i\|^2_R\big\}dt,
\end{align}
where $x_i\in \mathbb{R}^n$, $u_i\in \mathbb{R}^m$, $Q>0$ and $R>0$.

By Proposition \ref{lem2c}, it can be verified that A3)-A4) hold and (\ref{eq15c}) admits a solution $K_N>0$.
Furthermore, we have the following result.
\begin{proposition}
(i) (\ref{eq16c}) admits a unique $\rho$-stabilizing solution $P_N=0$; 

(ii) (\ref{eq17c}) admits a unique bounded solution
	$s_N(t)\equiv0$;

(iii) $ \mathbb{E}[\check{x}_i(t)]\equiv \bar{x}_0.$
\end{proposition}
\emph{Proof.} For the model (\ref{eq61})-(\ref{eq62}), (\ref{eq16c}) degenerates to
$\rho P_N=-P_N \Upsilon_N^{-1}P_N.$
It can be verified that $P_N=0$ is a unique $\rho$-stabilizing solution. Note that (\ref{eq17c}) reduces to
$\dot{s}_N=\rho s_N$. This with ${s}_N\in C_{\rho/2}([0,\infty),\mathbb{R}^n)$ implies	$s_N(t)\equiv0$ for any $t\geq 0$. Thus, 
we have $ \mathbb{E}[\check{x}_i(t)]\equiv \bar{x}_0.$
 \hfill$\Box$

For the model (\ref{eq61})-(\ref{eq62}),
 the decentralized strategies may be given as follows:
\begin{equation}\label{eq64a}
	\begin{split}
\check{u}_i(t)=&-\Upsilon_N^{-1}B^T{K}_N(\check{x}_i(t)-\bar{x}_0),  \ i=1,\cdots, N.
\end{split}
\end{equation}
Substituting (\ref{eq64a}) into (\ref{eq61}),
the closed-loop dynamics of agent $i$ can be written as
\begin{equation}\label{eq423}
	\begin{split}
		d\check{x}_i(t)=&-B\Upsilon_N^{-1}B^T{K}_N(\check{x}_i(t)-\bar{x}_0)dt\\&-D\Upsilon_N^{-1}B^T{K}_N(\check{x}_i(t)-\bar{x}_0)dw_i(t).
	\end{split}
  \end{equation}

 It can be shown that all the agents can achieve 
mean-square consensus.
  \begin{definition}\label{def1}	
	In a multiagent system, the agents  are said to reach
	the mean-square consensus if there exists a random variable
	$x^*$ such that
	$ \lim_{t\to \infty}\mathbb{E}\|x_i(t)-x^*\|^2=0$. 
	\end{definition}
\begin{theorem}\label{thm5.2}
For the model (\ref{eq61})-(\ref{eq62}), all the agents 
reach mean-square consensus. Specifically, under the strategy (\ref{eq64a}), there exist $c_1,c_2>0$ such that
$$ \mathbb{E}\|\check{x}^{(N)}(t)-\bar{x}_0\|^2\leq c_1 e^{-c_2 t}, \quad 1\leq i\leq N.$$

\end{theorem}

\emph{Proof.} See Appendix B. $\hfill \Box$

\section{Comparison and Discussion}

\subsection{Comparison with classical  
 mean-field games}

We now review the classical results of mean-field games for comparison to this work. Consider the large-population case of Problem (G) with 
 $\alpha_i^{(N)}=\frac{1}{N}$. By the mean-field (NCE) approach \cite{HCM07, WNZ19}, the following decentralized strategies are obtained
\begin{equation}\label{eq22}
{u}_i^*=-\Upsilon^{-1}\big[(B^TK+D^TKC)x_i+B^T\phi+D^TK\sigma\big],\quad
\end{equation}
where $ i=1,\cdots,N$, $\Upsilon=R+D^TKD$ and $K$ is the unique solution of the differential equation
\begin{equation}\label{eq23}
	\begin{split}
		\rho K&={A}^TK+KA+C^TKC-(B^TK+D^TKC)^T\cr&\times\Upsilon^{-1}(B^TK+D^TKC)+Q, \quad K(T)=0,
	\end{split}
  \end{equation}
   and $\phi$ is determined by the fixed-point equation
\begin{equation}\label{eq75}
\left\{  \begin{split}
\frac{d\phi}{dt}=&-[A-B\Upsilon^{-1}(B^TK+D^TKC)-\rho I]^T\phi\cr&-Kf+Q\Gamma(\bar{x}+\eta),\quad \phi(T)=0,\\ 
\frac{d\bar{x}}{dt}=&[A-B\Upsilon^{-1}(B^TK+D^TKC)]\bar{x}\\&-BR^{-1}B^T\phi+f, \ \bar{x}(0)= \bar{x}_0.
\end{split}
\right.
\end{equation}
Such set of decentralized strategies (\ref{eq22}) is further shown to be an $\varepsilon$-Nash equilibrium with respect to $\mathcal{U}_{c,i}$, i.e.,
  $$J_i(\check{u}_i,{u}_{-i}^*)\leq \inf_{u_i\in \mathcal{U}_{c,i}}J_i(u_i,\check{u}_{-i})+\varepsilon,$$
 where ${\mathcal U}_{c,i} =\big\{u_i, |u_i(t)\ \hbox{is adapted to}\ \sigma\{\bigcup_{i=1}^N{\mathcal F}_t^{i}\}\big\}$ and $\varepsilon=O(1/\sqrt{N})$.

Let $s(t)= \phi(t)-\Pi \bar{x}(t),$ where $\Pi\in \mathbb{R}^{n\times n}$. Then
\begin{align*}
\dot s=&\Pi\big[\big(A-B\Upsilon^{-1}(B^TK+D^TKC)\big)\bar{x}\cr&-B\Upsilon^{-1}B^T(\Pi\bar{x}+ s)+f\big]+\dot{ s}\cr
=&-\big(A-B\Upsilon^{-1}(B^TK+D^TKC)-\rho I\big)^T(\Pi\bar{x}+ s)\cr&-Kf+Q(\Gamma\bar{x}+\eta),
\end{align*}
where the second equation follows from (\ref{eq75}).
Comparing the terms in the equation above gives
\begin{align}\label{eq25}
\rho \Pi=& \Pi\big(A-B\Upsilon^{-1}(B^TK+D^TKC)\big)\cr&+\big(A-B\Upsilon^{-1}(B^TK+D^TKC)\big)^T\Pi\cr
&-\Pi B\Upsilon^{-1}B^T\Pi-Q\Gamma,\\
\label{eq26}
\rho s=&\dot{ s}+[A-B\Upsilon^{-1}(B^TK+B^T\Pi+D^TKC)]^T  s\cr&+(K+\Pi)f-Q\eta.
\end{align}

We have the following result.
\begin{proposition}
Assume that 
$\alpha_i^{(N)}=\frac{1}{N}$ for any $i=1,\cdots,N$. Then
   A2) holds for all sufficiently large $N$ and
$
    \|K_N\|+\|\Pi_N\|+\|s_N\|<\infty
  $
    if and only if (\ref{eq25})  admits a solution in $C([0,T],\mathbb{R}^{n\times n})$. Furthermore, we have
    \begin{equation}\label{eq31}
   \|K_N-K\|+\|\Pi_N-\Pi\|+\|s_N- s\|=o(1).
        \end{equation}
\end{proposition}
\emph{Proof.} If
(\ref{eq25}) admits a solution, then by the continuous dependence of the solution on the parameter $1/N$ (see e.g., \cite{F96}), there exists $N_0\geq 1$  such that (\ref{eq9b}) admits a solution for all $N\geq N_0$ (i.e., A2) holds), and (\ref{eq31}) is established.
This implies  $
     \|K_N\|+\|\Pi_N\|+\|s_N\|<\infty
  $ holds.
On other hand, if  Assumption A2) holds for all sufficiently large $N$ and
   $\|K_N\|+\|\Pi_N\|+\|s_N\|<\infty,$ then $\{\Pi_N,N\geq1\}$ are uniformly bounded. 
   Then there exists a subsequence $\{\Pi_{N_k},k\geq 1\}$ such that $\Pi_{N_k}$ converges to $\bar{\Pi}$ when $k\to \infty$. It can be verified that $\bar{\Pi}$ satisfie (\ref{eq25}). Thus, (\ref{eq25})  admits a solution. $\hfill \Box$




\begin{remark}
  The set of decentralized strategies (\ref{eq15a}) is an exact Nash equilibrium with respect to $\mathcal{U}_{d,i}, i=1,\cdots,N$. It is applicable for arbitrary number of agents.
In contrast, the set of decentralized strategies (\ref{eq22}) is an asymptotic Nash equilibrium with respect to $\mathcal{U}_{c,i}, i=1,\cdots,N$. It is only applicable for the large population case. However, the control gains of  
 (\ref{eq15a}) and (\ref{eq22}) coincide for the infinite population case.
\end{remark}

\subsection
{Computation of $\mathbb{E}[\check{x}_i]$ using consensus}

The decentralized strategies (\ref{eq15a}) actually involves coupling between agents due to the fact that $\mathbb{E}[\check{x}_i]$ satisfies (\ref{eq19}) which requires the averaged initial condition $\bar{x}_0$. For the infinite population case, the classical method is to compute $\bar{x}_0=\mathbb{E}[x_i(0)]$ by applying the statistical distribution of $x_i(0)$
or Monte-Carlo simulation. For the finite population case, we may use the average consensus algorithm to obtain $x^{(N)}(0)$. 
 Note that the asymptotic behavior of $\mathbb{E}[\check{x}_i]$ is not affected by the initial $\bar{x}_0$. 

Suppose there exist local interactions among all agents.
 Let a graph $(V,\mathcal{E})$ be given, where $V=\{1,2,\cdots,N\}$ is the set of vertices, and $\mathcal{E}=V\times V$ is the set of edges. Denote the set of neighbors of agent $i$ by $\mathcal{N}_i=\{j\in V|(i,j)\in \mathcal{E}\}$.
Given $x_i(0),i=1,\cdots,N$, we may utilize the average consensus algorithm to obtain $x^{(N)}(0)$: $$y_i(k+1)=y_i(k)+\sum_{j\in \mathcal{N}_i}l_{ij}(y_j(k)-y_i(k)), \ y_i(0)=x_i(0), $$
where $i=1,\cdots,N,$ and $L=(l_{ij})$ is the corresponding Laplacian matrix.
If $(V,\mathcal{E})$ is connected or has a spanning tree, then $y_i(k)$ will converge to $x^{(N)}(0)$, as $k\to \infty$. See e.g. \cite{KLM03, OFM07} for more details.
Applying a belief propagation (BP)-like distributed algorithm, the consensus will be reached with a fast convergence rate \cite{ZCF19}.



\section{ Finite-population Mean-field LQG Teams}\label{sec6}

In this section, we study the mean-field LQG social control problem with a finite number of agents. Both finite-horizon and infinite-horizon problems will be discussed.

\subsection{The finite-horizon problem}

We first consider the finite-horizon social problem.

{ \textbf{(S$^{\prime}$)}: Minimize social cost $J_{\rm soc, T}(u)=\sum_{i=1}^N\alpha_i^{(N)} J_{i,\rm T}$,} where
$$\left\{\begin{aligned}
J_{i,\rm T}(u)=&
\mathbb{E}\int_0^{T}e^{-\rho t}\Big[\big\|x_i-\Gamma x^{(\alpha)}-\eta\big\|^2_{Q}+\|u_i\|^2_R\Big]dt,\cr
 d{x}_i=&(A{x}_i+B{u}_i
		+ f)dt+
		(Cx_i+Du_i+\sigma)  dw_i.
\end{aligned}\right. $$
%
Denote
$$
\left\{
\begin{array}{l}
Q_{\Gamma}\stackrel{\Delta}{=}\Gamma^TQ+Q\Gamma-\Gamma^TQ\Gamma,\\
\bar{\eta}\stackrel{\Delta}{=}Q\eta-\Gamma^T Q\eta.
\end{array}\right.
$$
\begin{theorem}\label{thm1b}
{The problem (S$^{\prime}$) has a set of decentralized social optimal strategies $\hat{u}_i\in \mathcal{U}_{d,i}, i=1,\cdots,N$, if and only if the following equation system admits a set of solutions $(\hat{x}_i,\hat{\lambda}_i, \hat{\beta}_{i}^{j},i,j=1,\cdots,N)$:
 \begin{equation}
\label{eq63}
\left\{ \begin{aligned}
d\hat{x}_i=&(A\hat{x}_i+B\hat{u}_i+f)dt
  +(C\hat{x}_i+D\hat{u}_i+\sigma) dw_i,\cr
d\hat{\lambda}_i=&-\big[(A-\rho I)^T\hat{\lambda}_i+C^T\hat{\beta}_i^i+Q\hat{x}_i-Q_{\Gamma} \hat{x}^{(\alpha)}-\bar{\eta}\big]dt\cr&+ \sum_{j=1}^N
  \hat{\beta}_{i}^j dw_j, \\
\ x_i(0) &=x_{i0},\quad \hat{\lambda}_i(T)=0, i=1,\cdots,N,
\end{aligned}\right.
 \end{equation}}
with
\begin{equation}\label{eq7c}
	\begin{split}		 \hat{u}_i=&-R^{-1}\big(B^T\mathbb{E}[\hat{\lambda}_i|\mathcal{F}_t^i]+D^T\mathbb{E}[\hat{\beta}^i_i|\mathcal{F}_t^i]\big).
	\end{split}
  \end{equation}
\end{theorem}
\emph{Proof.} See Appendix C. \hfill{$\Box$}

For simplicity, consider 
the case $\alpha_i^{(N)}=\frac{1}{N}$ later.
 Recall  $$\mathbb{E}[\hat{x}^{(N)}|\mathcal{F}_t^i]=\frac{1}{N}\hat{x}_i+\frac{N-1}{N}\mathbb{E}[\hat{x}_i].$$
It follows from (\ref{eq63}) that
\begin{equation}\label{eq65}
\left\{
\begin{aligned}
d\mathbb{E}[\hat{x}_i]=&(A\mathbb{E}[\hat{x}_i]+B \mathbb{E}[\hat{u}_i]+f)dt, \mathbb{E}[\hat{x}_i(0)]=\bar{x}_{0},\cr
    d\mathbb{E}[\hat{\lambda}_i|\mathcal{F}_t^i]=& - \big[(A-\rho I)^T\mathbb{E}[\hat{\lambda}_i|\mathcal{F}_t^i]+\big(Q-\frac{1}{N}Q_{\Gamma}\big)\hat{x}_i\cr&-\frac{N-1}{N}Q_{\Gamma}\mathbb{E}[\hat{x}_i]-\bar{\eta} \big]dt-C^T\mathbb{E}[\hat{\beta}^i_i|\mathcal{F}_t^i]dt\cr&+\mathbb{E}[\hat{\beta}^i_i|\mathcal{F}_t^i]dw_i, \ \mathbb{E}[\hat{\lambda}_i(T)|\mathcal{F}_T^i]=0.
\end{aligned}\right.
\end{equation}
  Let
 $\mathbb{E}[\hat{\lambda}_i(t)|\mathcal{F}_t^i]=\hat{K}_N(t)\hat{x}_i(t)+\hat{\Pi}_N(t)\mathbb{E}[\hat{x}_i(t)]+\hat{s}_N(t),$
     where $\hat{K}_N(t),\hat{\Pi}_N(t)\in \mathbb{R}^{n\times n}$ and $\hat{s}_N(t)\in \mathbb{R}^n$.
 By It\^{o}'s formula, we obtain
\begin{align}\label{eq66}
d\mathbb{E}[\hat{\lambda}_i|\mathcal{F}_t^i]\!=&\dot{\hat{K}}_N\hat{x}_idt+\hat{K}_N\big[(A\hat{x}_i+B\hat{u}_i+f)dt\cr
+&(C\hat{x}_i+D\hat{u}_i+\sigma) dw_i\big]+\dot{\hat{\Pi}}_N\mathbb{E}[\hat{x}_i]dt
\cr\!+&\hat{\Pi}_N\big[A\mathbb{E}[\hat{x}_i]+B \mathbb{E}[\hat{u}_i]+f\big]dt+\dot{\hat{s}}_Ndt.
\end{align}
Comparing this with (\ref{eq65}), it follows that $\mathbb{E}[\hat{\beta}^i_i|\mathcal{F}_t^i]=\hat{K}_N
(C\hat{x}_i+D\hat{u}_i+\sigma)$.
By (\ref{eq7c}), we have
\begin{equation}\label{eq67}
	\begin{split}
		\hat{u}_i=&-\hat{\Upsilon}_N^{-1}\big[(B^T\hat{K}_N+D^T\hat{K}_NC)\hat{x}_i\\&+B^T\hat{\Pi}_N\mathbb{E}[\hat{x}_i]+B^T\hat{s}_N+D^T\hat{K}_N\sigma\big],
	\end{split}
  \end{equation}
where $\hat{\Upsilon}_N\stackrel{\Delta}{=}R+D^T\hat{K}_ND.$
Applying (\ref{eq67}) to (\ref{eq63}), it follows that
\begin{align}\label{eq68}
\rho \hat{K}_N\!=&\dot{\hat{K}}_N+{A}^T\hat{K}_N+\hat{K}_NA+C^T\hat{K}_NC\cr-&
(B^T\hat{K}_N+D^T\hat{K}_NC)^T\hat{\Upsilon}_N^{-1}(B^T\hat{K}_N+D^T\hat{K}_NC)\cr
+& Q-\frac{1}{N}Q_{\Gamma}, \quad \hat{K}_N(T)=0,\\
\label{eq69}
\rho \hat{\Pi}_N\!=&
 \dot{\hat{\Pi}}_N+{A}^T\hat{\Pi}_N+\hat{\Pi}_N{A}-\hat{\Pi}_N B \hat{\Upsilon}_N^{-1}B^T\hat{\Pi}_N\cr-&\hat{\Pi}_N B \hat{\Upsilon}_N^{-1}(B^T\hat{K}_N+D^T\hat{K}_NC)-\frac{N-1}{N}Q_{\Gamma}\cr
-&(B^T\hat{K}_N+D^T\hat{K}_NC)^T\hat{\Upsilon}_N^{-1}B^T\hat{\Pi}_N, \ \hat{\Pi}_N(T)=0,\\
\label{eq70}
\rho \hat{s}_N\!=&\dot{\hat{s}}_N\!+\!\big[{A}\!-\!B\hat{\Upsilon}_N^{-1}\big( B^T(\hat{K}_N+\hat{\Pi}_N)+D^T\hat{K}_NC\big)\big]\hat{s}_N\cr
+&\big[C\!-D\hat{\Upsilon}_N^{-1}\big( B^T(\hat{K}_N+\hat{\Pi}_N)+D^T\hat{K}_NC\big)\big]^T\hat{K}_N\sigma,\cr+&(\hat{K}_N+\hat{\Pi}_N)f-\bar{\eta},\
\hat{s}_N(T)=0.
\end{align}

%

We have the following result.

\begin{proposition}\label{prop4}
  Equations (\ref{eq68})-(\ref{eq69}) admit solutions in $C([0,T], \mathbb{R}^{n\times n})$.
\end{proposition}

\emph{Proof.} See Appendix C. \hfill{$\Box$}

Note that by Proposition \ref{prop4}, (\ref{eq68})-(\ref{eq70}) admit a solution, respectively.
Then we obtain the decentralized strategy (\ref{eq67}),
where $\hat{K}_N, \hat{\Pi}_N, \hat{s}_N$ are given by (\ref{eq68})-(\ref{eq70}).  
The following theorem gives the performance of the proposed decentralized strategy above.
\begin{theorem} \label{thm6.2}
	Let A1) hold.
Then for Problem (S$^{\prime}$), the set of control strategies
	$\{\hat{u}_1,\cdots,\hat{u}_N\}$ given by (\ref{eq67}) is  a decentralized social optimal solution, 
  and the corresponding social cost is given by
  \begin{align}\label{eq16e}
&J_{\rm soc,T}(\hat{u})=\sum_{i=1}^N\mathbb{E}\Big\{\big\|x_{i0}-\mathbb{E}[x_i(0)]\big\|^2_{\hat{K}_N(0)}\cr
&+\big\|\mathbb{E}[x_i(0)]\big\|^2_{\hat{P}_N(0)} +2s^T_N(0)\mathbb{E}[x_i(0)]\Big\}+Nq_T^N,
  \end{align}	
  where
  \begin{align}\label{eq17a}
  q_T^N=&\int_0^{T}e^{-\rho t}\big[\|\sigma \|^2_{\hat{K}_N}+\| \sigma\|^2_{\hat{P}_N} +\|\eta\|^2_Q \cr&-\|B^Ts_N+D^T\hat{K}_N\sigma\|^2_{{\hat{\Upsilon}}_N^{-1}}+2s^T_Nf\big]dt.
  \end{align}
\end{theorem}
\emph{Proof.} See Appendix C. \hfill{$\Box$}

\begin{remark}
The work \cite {WZZ20} investigated mean-field social control for the case $C=D=0$. When the classical mean-field social controller $\{u_i^*=-\Upsilon^{-1}B^T(Kx_i+(P-K)\bar{x}+s),i=1,\cdots,N\}$ 
 is applied, the corresponding social cost is given by
 \begin{align*}
J_{\rm soc, T}(\hat{u})=&\sum_{i=1}^N\mathbb{E}\Big\{\big\|x_{i0}-x^{(N)}(0)\big\|^2_{K(0)}+\big\|x^{(N)}(0)\big\|^2_{P(0)}\cr&
  +2s^T(0)x^{(N)}(0)\Big\}+Nq_T+N\epsilon_T,
  \end{align*}
where
\begin{align*}
  q_T=&\int_0^{T}e^{-\rho t}\big[\|\sigma \|^2_{K}+\| \sigma\|^2_{P} -\|B^Ts\|^2_{R^{-1}}+2s^Tf\big]dt,\\
\epsilon_T=&\mathbb{E}\int_0^{T}\!e^{-\rho t}\big\|B^T[P-K][x^{(N)}\!-\!\bar{x}]\big\|^2_{R^{-1}}dt.
  \end{align*}
  Here, $K$ satisfies (\ref{eq23}) and $P$ satisfies
  \begin{equation*}
  \begin{aligned}
\rho {P}=  &\dot{P}+{A}^TP+PA-P B
  R^{-1}B^TP\cr&+(I-\Gamma)^TQ(I-\Gamma),\quad  P(T)=0.
  \end{aligned}
  \end{equation*}
  Compared to (\ref{eq16e}), residual term $\epsilon_T$ vanishes as $N\to\infty$.
\end{remark}

\subsection{The Infinite-Horizon Problem}
Based on the discussion above, we may design the following decentralized strategy:
\begin{equation}\label{eq75b}
	\begin{split}		 \hat{u}_i=&\!-\!\hat{\Upsilon}_N^{-1}\big[(B^T\hat{K}_N\!+\!D^T\hat{K}_NC)\hat{x}_i\!+\!B^T(\hat{P}_N-\hat{K}_N)\mathbb{E}[\hat{x}_i]\cr&+B^T\hat{s}_N+D^T\hat{K}_N\sigma\big], \quad i=1,\cdots, N,
	\end{split}
\end{equation}
where $\hat{\Upsilon}_N=R+D^T\hat{K}_ND$ and $\hat{K}_N, \hat{P}_N, \hat{s}_N$ are given by
\begin{align}\label{eq76}
\rho \hat{K}_N\!\!=&{A}^T\hat{K}_N+\hat{K}_NA-(B^T\hat{K}_N+D^T\hat{K}_NC)^T\hat{\Upsilon}_N^{-1}\cr&\times(B^T\hat{K}_N+D^T\hat{K}_NC)+C^T\hat{K}_NC+ Q-\frac{1}{N}Q_{\Gamma},\\
\label{eq77}
  \rho \hat{P}_N\!\!=&{A}^T\hat{P}_N+\hat{P}_N{A}-(B^T\hat{P}_N+D^T\hat{K}_NC)^T \hat{\Upsilon}_N^{-1}\cr&\times(B^T\hat{P}_N+D^T\hat{K}_NC)+C^T\hat{K}_NC+Q-Q_{\Gamma},\\
\label{eq78}
\rho \hat{s}_N\!\!=&\dot{\hat{s}}_N+\big[{A}-B\hat{\Upsilon}_N^{-1} (B^T\hat{P}_N^T+D^T\hat{K}_NC\big)\big]\hat{s}_N+\hat{P}_Nf\cr
&-\bar{\eta}+\big[C-D\hat{\Upsilon}_N^{-1}\big( B^T\hat{P}_N^T+D^T\hat{K}_NC\big)\big]^T\hat{K}_N\sigma.
\end{align}

For further analysis, we introduce the assumption: 

\textbf{A5)}  (\ref{eq77}) admits a unique $\rho$-stabilizing solution.

Note $$Q-\frac{1}{N}Q_{\Gamma}=\frac{N-1}{N}Q+\frac{1}{N}(I-\Gamma)^TQ(I-\Gamma).$$ From A3) and \cite[Theorem 4.1]{ZZC08}, we obtain that
$\big[A-\frac{\rho}{2}I, C,\sqrt{Q-\frac{1}{N}Q_{\Gamma}}\big]$ is exactly detectable and (\ref{eq76}) admits a unique $\rho$-stabilizing solution. Define
$$\widehat{\mathcal{M}}_{\Gamma}{=}\left[\begin{array}{cc}
		\hat{A}_N-\frac{\rho}{2}I&B \hat{\Upsilon}_N^{-1}B^T\\
		C^T_N+\big(Q-Q_{\Gamma}\big)  &-\hat{A}_N^{T}+\frac{\rho}{2}I
	\end{array}
	\right].$$

The following proposition provides some conditions to ensure that A5) holds.
\begin{proposition}\label{prop6.2}
  Assume that A3) holds. Then  A5) holds if and only if $\widehat{\mathcal{M}}_{\Gamma}$ has no eigenvalues on the imaginary axis.
  Particularly, if (i) $[A-\frac{\rho}{2}I, C,\sqrt{Q(I-\Gamma)}]$ is exactly detectable or (ii) $\Gamma=I$, $C=0$ and $A-\frac{\rho}{2}I$ has no eigenvalues on the imaginary axis, then A5) holds.
\end{proposition}

\emph{Proof.} Since A3) holds, by \cite{M77}, A5) holds if and only if $\widehat{\mathcal{M}}_{\Gamma}$ has no eigenvalues on the imaginary axis.
 Particularly, if $[A-\frac{\rho}{2}I, C,\sqrt{Q(I-\Gamma)}]$ is exactly detectable, then from A3) and \cite[Theorem 4.1]{ZZC08}, (\ref{eq77}) admits a unique $\rho$-stabilizing solution.
 Note $Q-Q_{\Gamma}=0$ for $\Gamma=I$. If $C=0$, and $A-\frac{\rho}{2}I$ has no eigenvalues on the imaginary axis, then $\widehat{\mathcal{M}}_{\Gamma}$ has no eigenvalues on the imaginary axis, which further implies that A5) holds. \hfill{$\Box$}

%

\begin{remark}
  For the  case $\Gamma=I$ and $\eta=0$, we have $\big(I-\frac{1}{N}\Gamma\big)^TQ(I-\Gamma)=Q-Q_{\Gamma}=0$, which gives that (\ref{eq16c}) and (\ref{eq77}) have the same solutions. On other hand,
$$(I-\frac{1}{N})^TQ(I-\frac{1}{N})=\frac{(N-1)^2}{N^2}Q\not =\frac{N-1}{N}Q=Q-\frac{1}{N}Q_{\Gamma}.$$ This implies the solutions to (\ref{eq15c}) and (\ref{eq76}) are different somewhat.
 Thus, the social solution and the game solution are slightly different for the finite-population consensus problem. However, the two solutions coincide for the infinite population case.
\end{remark}

 Similar to Theorem \ref{thm6.2}, we obtain the following result.
\begin{theorem} 
	Let A1), A3) and A5) hold.
Then for Problem (S), the set of control laws
	$\{\hat{u}_1,\cdots,\hat{u}_N\}$ given by (\ref{eq75b}) is a decentralized social optimal solution.
\end{theorem}

\subsection{Extension to the case with state coupling}

{Consider the case that agents $ i=1,\cdots, N$ evolves by
		\begin{eqnarray}\label{eq55}
			dx_i\!=\!(\!Ax_i\!+\!Bu_i\!+Gx^{(N)}+\!f)dt\!+(Cx_i+\sigma)dw_i,
		\end{eqnarray}
with the social cost
		\begin{align} \label{eq56}
				J_{\rm soc,T}(u)=\sum_{i=1}^N\mathbb{E}\int_0^Te^{-\rho t} (|x_i-\Gamma x^{(N)}-\eta|^2_Q+|u_i|^2_R)dt.
		\end{align}
}

{By a similar derivation for Theorem \ref{thm1b}, we obtain that the problem (\ref{eq55})-(\ref{eq56}) has a set of team-optimal strategies $\breve{u}_i\in \mathcal{U}_{d,i}, i=1,\cdots,N$ if and only if the FBSDEs admit a  solution $(\breve{x}_i,\breve{\lambda}_i, \breve{\beta}_{i}^{j},i,j=1,\cdots,N)$:
  	\begin{equation*}
  			\left\{\begin{aligned}
  		d\breve{x}_i=&(A\breve{x}_i+B\breve{u}_i+{G\breve{x}^{(N)}}+f)dt+(C\breve{x}_i+\sigma)dw_i, \\
  		d\breve{\lambda}_i=&	-\Big\{[A-\rho I]^T\breve{\lambda}_i+C^T\breve{\beta}_i^i+{G^T\breve{\lambda}^{(N)}}+Q\breve{x}_i\\&-Q_\Gamma \breve{x}^{(N)}-\bar{\eta}
  		\Big\}dt+\sum_{j=1}^N\breve{\beta}_i^jdw_j,\cr
 \breve{x}_i=&x_{i0},\    \breve{\lambda}_i(T)=0, \ i=1, \cdots,N,
  	\end{aligned}\right.
  	\end{equation*}
  	and furthermore the optimal control laws are given by $$	\breve{u}_i=-R^{-1}B^T\mathbb{E}[\breve{\lambda}_i|\mathcal{F}_t^i],\  i=1,\cdots,N.$$
Note that for $1\leq i\not = j\leq N$,
\vspace{-2mm}
$${\mathbb{E}[\breve{x}^{(N)}|\mathcal{F}_t^i]=\frac{1}{N}\mathbb{E}[\breve{x}_i|\mathcal{F}_t^i]+\frac{N-1}{N}\mathbb{E}[\breve{x}_j|\mathcal{F}_t^i],}$$
where ${\mathcal F}_t^{i}=\sigma(x_{i}(0),w_i(s),0\leq s\leq t)$. We have
{ \small \begin{equation*}
	\left\{\begin{aligned}
		d\mathbb{E}[\breve{x}_i|\mathcal{F}_t^i]=&\Big\{\big(A+\frac{G}{N}\big)\mathbb{E}[\breve{x}_i|\mathcal{F}_t^i]-BR^{-1}
		B^T  \mathbb{E}[\breve{\lambda}_i|\mathcal{F}_t^i]+f\\
+&\frac{N-1}{N}G\mathbb{E}[\breve{x}_j|\mathcal{F}_t^i]\Big\}dt+ (C\mathbb{E}[\breve{x}_i|\mathcal{F}_t^i]+\sigma) dw_i, \\	d\mathbb{E}[\breve{x}_j|\mathcal{F}_t^i]=&\Big\{\big(A+\frac{N-1}{N}G\big)\mathbb{E}[\breve{x}_j|\mathcal{F}_t^i]+
\frac{G}{N}\mathbb{E}[\breve{x}_i|\mathcal{F}_t^i]\\
-&BR^{-1}
		B^T  \mathbb{E}[\breve{\lambda}_j|\mathcal{F}_t^i]+f\Big\}dt,\\
		d\mathbb{E}[\breve{\lambda}_i|\mathcal{F}^i_t]=&	-\Big\{\big[A-\rho I+\frac{G}{N}\big]^T\mathbb{E}[\breve{\lambda}_i|\mathcal{F}_t^i]+C^T\mathbb{E}[\breve{\beta}_i^i|\mathcal{F}_t^i]
\\+&\frac{N-1}{N}G^T\mathbb{E}[\breve{\lambda}_j|\mathcal{F}_t^i]+\big(Q-\frac{Q_\Gamma}{N}\big)\mathbb{E}[\breve{x}_i|\mathcal{F}_t^i]\\
-&\frac{N-1}{N}Q_\Gamma\mathbb{E}[\breve{x}_j|\mathcal{F}_t^i]-\bar{\eta}\Big\}dt+\mathbb{E}[\breve{\beta}_i^i|\mathcal{F}_t^i]dw_i,\\
		d\mathbb{E}[\breve{\lambda}_j|\mathcal{F}_t^i]=&	-\Big\{\big[A-\rho I+\frac{N-1}{N}G\big]^T\mathbb{E}[\breve{\lambda}_j|\mathcal{F}_t^i]\\
+&\frac{G^T}{N}\mathbb{E}[\breve{\lambda}_i|\mathcal{F}_t^i]
+C^T\mathbb{E}[\breve{\beta}_j^j|\mathcal{F}_t^i]-\frac{Q_\Gamma}{N}\mathbb{E}[\breve{x}_i|\mathcal{F}_t^i]\\
+&\big(Q\!-\!\frac{N-1}{N}Q_\Gamma\big)\mathbb{E}[\breve{x}_j|\mathcal{F}_t^i]\!-\!\bar{\eta} \Big\}dt\!+\!\mathbb{E}[\breve{\beta}_j^i|\mathcal{F}_t^i]dw_i,\\
		\mathbb{E}[\breve{x}_i(0)|\mathcal{F}^i_t]&=x_{i0},\mathbb{E}[\breve{x}_j|\mathcal{F}^i_t]=\bar{x}_{0}, \mathbb{E}[\breve{\lambda}_i(T)|\mathcal{F}_t^i]=0,
\\&\mathbb{E}[\breve{\lambda}_j(T)|\mathcal{F}_t^i]=0, \quad  1\leq i\not = j\leq N.
	\end{aligned}\right.
\end{equation*}}
Let $\breve{\lambda}_i=\breve{K}_N\breve{x}_i+\breve{\Pi}_N\breve{x}^{(N)}+\breve{s}_N$. 
 We have 		 {\small 	
 \begin{align*}	\mathbb{E}[\breve{\lambda}_i|\mathcal{F}_t^i]&=(\breve{K}_N+\frac{\breve{\Pi}_N}{N})\mathbb{E}[\breve{x}_i|\mathcal{F}_t^i]+\frac{N-1}{N}\breve{\Pi}_N \mathbb{E}[\breve{x}_j|\mathcal{F}_t^i]+\breve{s}_N,\\			 			  			\mathbb{E}[\breve{\lambda}_j|\mathcal{F}_t^i]&=\frac{\breve{\Pi}_N}{N}\mathbb{E}[\breve{x}_i|\mathcal{F}_t^i]+(\breve{K}_N+\frac{N-1}{N}\breve{\Pi}_N) \mathbb{E}[\breve{x}_j|\mathcal{F}_t^i]+\breve{s}_N.
 \end{align*}}
By applying
 the four-step scheme \cite{MY99}, we obtain  the following social optimal control
$$\begin{aligned}
  \breve{u}_i=&-R^{-1}B^T\big[(\breve{K}_N+\frac{\breve{\Pi}_N}{N})\mathbb{E}[\breve{x}_i|\mathcal{F}_t^i\big]\cr
  &+\frac{N-1}{N}\breve{\Pi}_N \mathbb{E}[\breve{x}_j|\mathcal{F}_t^i]+\breve{s}_N],
  \end{aligned}$$
 where
\begin{align*}
  		\!\rho \breve{K}_N\!=&\dot{\breve{K}}_N+A^T\breve{K}_N+\breve{K}_NA+Q-\breve{K}_NBR^{-1}B^T\breve{K}_N\\
 &+C^T\big(\breve{K}_N+\frac{\breve{\Pi}_N}{N}\big)C, \ \breve{K}_N(T)=0,	\\
 		\!\rho \breve{\Pi}_N\!=&\dot{\breve{\Pi}}_N+A^T\breve{\Pi}_N+\breve{\Pi}_N A-\breve{K}_NB^TR^{-1}B\breve{\Pi}_N\\&-\breve{\Pi}_N BR^{-1}B^T(\breve{K}_N+\breve{\Pi}_N)+G^T(\breve{K}_N+\breve{\Pi}_N)\\&+(\breve{K}_N+\breve{\Pi}_N)G-Q_\Gamma,\ 
 \breve{\Pi}_N(T)=0,
 \end{align*}
 \begin{align*}
 		\!\rho \breve{s}_N\!=&\dot{\breve{s}}_N+(A+G-B^TR^{-1}B(\breve{K}_N
 +\breve{\Pi}_N))^T\breve{s}_N-\bar{\eta}\cr&+(\breve{K}_N
 +\breve{\Pi}_N)f+C^T(\breve{K}_N+\frac{\breve{\Pi}_N}{N})\sigma, \breve{s}_N(T)=0.
 	\end{align*}
}

\section{Numerical Examples}
In this section, we give two numerical examples to illustrate the properties of proposed decentralized strategies.

Consider Problem (G) for $6$ agents with single-integrator dynamics and additive noise (i.e. $A=0$ and $C=D=0$). Take the parameters as $B=Q=R=\Gamma= 1, f(t)=\eta(t)=0,
\sigma=0.1, \rho=0.2$, and $\alpha_i^{(N)}=\frac{1}{N}$.
The initial states of $6$ agents are taken independently from a normal distribution $N(5,1)$.  Note that $B\not=0$, and $Q>0$. Then
A1) and A3) hold. By Proposition \ref{lem2c}, A4) holds.

\begin{figure}[H]
	\centering
	\includegraphics[width=0.9\linewidth]{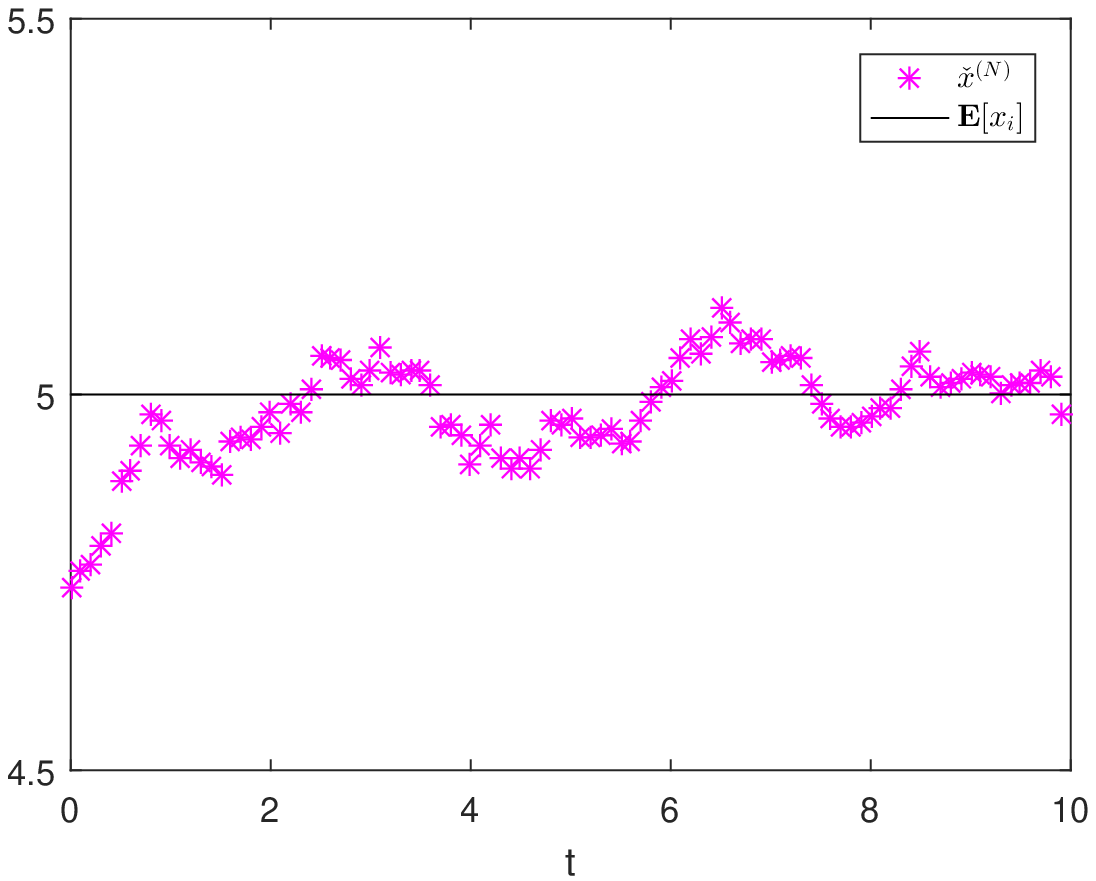}
	\caption{Curves of $\check{x}^{(N)}$ and $\mathbb{E}[\check{x}_i]$.}
	\label{consistency}
\end{figure}

Under the strategy (\ref{eq14c}),
the trajectories of $\mathbb{E}[\check{x}_i]$ and $\check{x}^{(N)}$ are shown in Fig. \ref{consistency}.
It can be seen that $\mathbb{E}[\check{x}_i]$ and $\check{x}^{(N)}$ do not coincide well, but $\mathbb{E}[{x}_i]$ attains the mean of $\check{x}^{(N)}$. This is different from classical mean-field games, where 
$\mathbb{E}[\check{x}_i]$ and $\check{x}^{(N)}$ coincide as agent number is large.

\begin{figure}[H]
	\centering
	\includegraphics[width=1.0\linewidth]{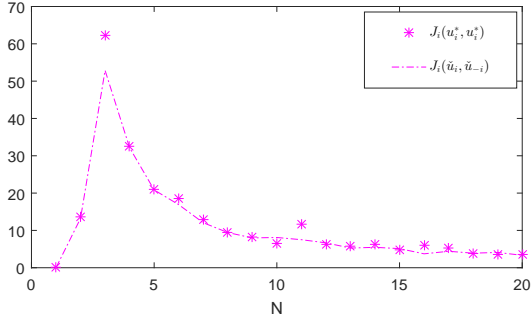}
	\caption{Performance comparison of proposed Nash strategy and classical  mean-field strategy.}
	\label{comparison}
\end{figure}

It can be seen from  Fig. \ref{comparison} that  the performance comparison of the proposed Nash equilibrium strategy and the classical  mean-field strategy.
When $N$ is not very large, it gives superior performance using the proposed Nash equilibrium (\ref{eq14c}) than the $\varepsilon$-Nash equilibrium given by the classical mean-field strategy (\ref{eq22}). The two performances merge when $N\to \infty$. 

\begin{figure}[H]
	\centering
	\includegraphics[width=0.9\linewidth]{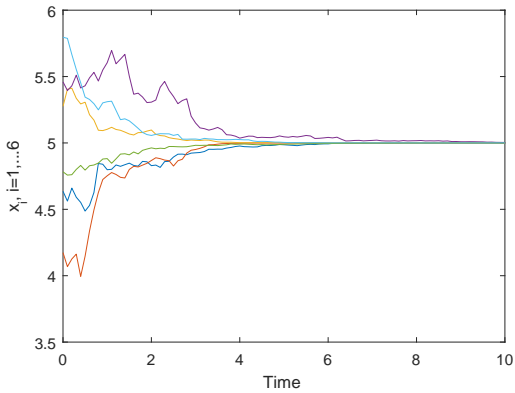}
	\caption{Curves of $\hat{x}_i$, $i=1,\cdots,6$.}
	\label{states2}
\end{figure}

Consider Problem (S) for $6$ single-integrator agents with multiplicative noise (i.e., $A=C=0$).
Take the parameters as $B=Q=R=D=\Gamma= 1, f(t)=\sigma(t)=\eta(t)=0, \rho=0.2$. 
The initial states of $6$ agents are taken independently from a normal distribution $N(5,1)$.  Note that $B=D\not=0$, and $Q>0$. Then
A1) and A3) hold. By Proposition \ref{prop6.2}, A5) holds.

\begin{figure}[H]
	\centering
	\includegraphics[width=0.9\linewidth]{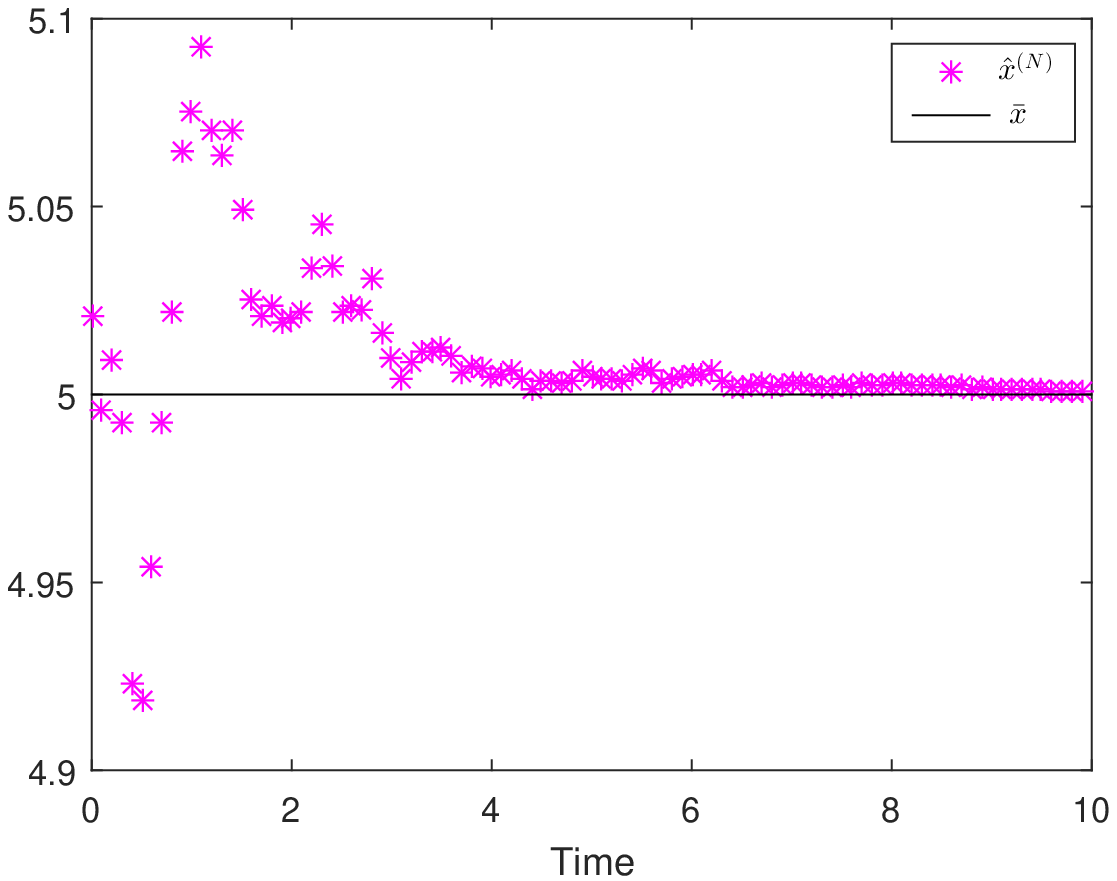}
	\caption{Curves of $\check{x}^{(N)}$ and $\mathbb{E}[\check{x}_i]$.}
	\label{consistency2}
\end{figure}

\begin{figure}[H]
	\centering
	\includegraphics[width=0.9\linewidth]{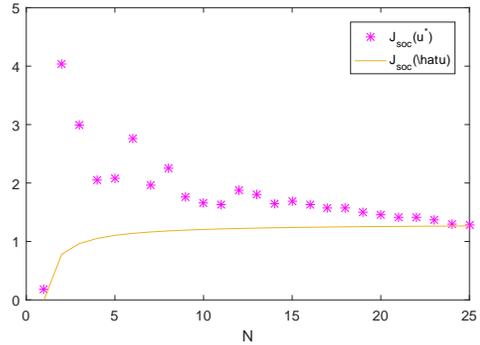}
	\caption{Performance comparison of proposed social strategy and classical mean-field controller.}
	\label{comparison2}
\end{figure}

Under the proposed control strategy (\ref{eq75b}),
the trajectories of $\hat{x}_i,i=1,\cdots,6$ are shown in Fig. \ref{states2}.
It can be seen that under multiplicative noise, all the agents can achieve consensus strictly, which verifies the result of Theorem \ref{thm5.2}.
The trajectories of $\mathbb{E}[\hat{x}_i]$ and $\hat{x}^{(N)}$ are shown in Fig. \ref{consistency2}.
It can be seen that $\hat{x}^{(N)}$ and $\mathbb{E}[\hat{x}_i]$ coincide well when the time is sufficiently long. 
This is different from the case with additive noise, which is shown in  Fig. \ref{consistency}.
Fig. \ref{comparison2} shows  the performance comparison of the proposed social strategy and the classical  mean-field social controller.
When $N$ is not very large, it gives superior performance using the proposed social strategy (\ref{eq67}) than 
 the classical  mean-field controller (\ref{eq22}). The two performances merge when $N$ is sufficiently large.

\section{Conclusion}

In this paper, we investigated mean-field games and teams for a finite number of agents. For finite-horizon problems, we designed  decentralized
strategies in terms of two differential Riccati equations by decoupling non-standard
FBSDEs. The proposed  decentralized strategies were further shown to be a Nash equilibrium and a social optimal solution, respectively. For infinite-horizon problems, we gave some criteria for the solvability of algebraic Riccati equations arising from
consensus. For future investigation, it would be interesting to
generalize the results to more complicated situations, such as mixed games with a major player or leader-follower games, and nonlinear games with a finite number of agents.


\appendix
\section{Proofs for Section \ref{sec3a}}
\def\theequation{A.\arabic{equation}}
\setcounter{equation}{0}

\emph{Proof of Theorem \ref{thm1}.}
 (i) Suppose $\{\check{u}_i,i=1,\cdots,N\}$ is a set of decentralized Nash equilibrium strategies of Problem (G$^{\prime}$), and $\{\check{x}_i,i=1,\cdots,N\}$ are the corresponding states of agents, i.e., they satisfy (\ref{eq6a}).
Let $\{\check{\lambda}_i, \check{\beta}_i^j,i,j=1,\cdots,N\}$ be a set of adapted solutions to the second equation of (\ref{eq6a}).
 For any $u_i\in \mathcal{U}_{d,i} 
 $ and $\theta\in \mathbb{R}\ (\theta \not= 0)$, let $v_i^{\theta}=\check{u}_i+\theta v_i$. Denote by $x_i^{\theta}$ the solution of the following perturbed state equation
$$ \begin{aligned}
  dx_i^{\theta}=&\big[Ax_i^{\theta}+B(\check{u}_i+\theta v_i)+f\big]dt
  +\big[Cx_i^{\theta}+D(\check{u}_i\cr&+\theta v_i)+\sigma\big] dw_i,\quad
 x_i^{\theta}(0)=x_{i0},\ i=1,2,\cdots,N.
 \end{aligned}$$ Let $y_i=(x_i^{\theta}-\check{x}_i)/\theta$. 
It can be verified that
 $y_i$ satisfies
\begin{equation}\label{eq4aa}
  \begin{aligned}
   &dy_i=(A y_i+Bv_i)dt+(Cy_i+Dv_i)dw_i,\ y_i(0)=0.
  \end{aligned}
\end{equation}
Then by It\^{o}'s formula, for any $i=1,\cdots,N$,
\begin{align}\label{eq7}
0
= &\mathbb{E}[\langle \check{\lambda}_i(T), e^{-\rho T}y_i(T)\rangle-\langle \check{\lambda}_i(0),y_i(0)\rangle]
\cr
=& \mathbb{E}\int_0^T e^{-\rho t}  \Big[\Big\langle -\big(I-\alpha_i^{(N)}\Gamma\big)^T Q(\check{x}_i-\Gamma \check{x}^{(\alpha)}-\eta),y_i\Big\rangle\cr&+\langle B^T\check{\lambda}_i+D^T\beta_i^i,v_i\rangle\Big] dt.
\end{align}
We have
\begin{equation}\label{eq5a}
\begin{aligned}
 &{J}_{i,\rm T}(u_i^{\theta},\check{u}_{-i})-{J}_{i, \rm T}(\check{u}_i, \check{u}_{-i})=2\theta I_1+{\theta^2}I_2
\end{aligned}
\end{equation}
where $\check{u}_{-i}=(\check{u}_1,\cdots,\check{u}_{i-1},\check{u}_{i+1},\cdots,\check{u}_N)$, and
\begin{align*}
I_1\stackrel{\Delta}{=}&\mathbb{E}\int_0^T  e^{-\rho t}\big[\big\langle Q\big(\check{x}_i-(\Gamma\check{x}^{(\alpha)}+\eta)\big),(I-\alpha_i^{(N)}\Gamma) y_i\big\rangle \cr&+
\langle R \check{u}_i,v_i\rangle \big]dt,\cr
I_2\stackrel{\Delta}{=}&\mathbb{E}\int_0^T e^{-\rho t}\big[\big\|(I-\alpha_i^{(N)}\Gamma)y_i
   \big\|^2_{Q}
+\|v_i\|^2_{R}\big]dt.
\end{align*}
From (\ref{eq7}), one can obtain that
\begin{align*}
 I_1=&\mathbb{E}\int_0^T  e^{-\rho t}\big[\big\langle Q\big(\check{x}_i-(\Gamma\check{x}^{(\alpha)}+\eta)\big),(I-\alpha_i^{(N)}\Gamma) y_i\big\rangle\cr& +
\langle R \check{u}_i,v_i\rangle \big]dt,+ \mathbb{E}\int_0^T e^{-\rho t}  \Big[\Big\langle -\Big(I-\alpha_i^{(N)}\Gamma\Big)^T Q\cr&\times(\check{x}_i-\Gamma \check{x}^{(\alpha)}-\eta),y_i\Big\rangle+\langle B^T\check{\lambda}_i+D^T\check{\beta}_i^i,v_i\rangle\Big] dt\cr
=&\mathbb{E}\int_0^T\Big\langle R\check{u}_i+B^T\check{\lambda}_i+D^T\check{\beta}_i^i,v_i\Big\rangle dt.
\end{align*}
Note that $\check{u}_i,v_i\in \mathcal{U}_{d,i}$. By the smoothing property of conditional mathematical expectation,
\begin{align}
  \label{eq10d}
 I_1=&\mathbb{E}\int_0^T\Big\langle R\check{u}_i+B^T\mathbb{E}[\check{\lambda}_i|\mathcal{F}_t^i]+D^T\mathbb{E}[\check{\beta}^i_i|\mathcal{F}_t^i],v_i\Big\rangle dt.
 \end{align}
 Since $Q \geq0$ and $R> 0$, we have
 $I_2\geq0$.
From (\ref{eq5a})-(\ref{eq10d}), the fact that $\{\check{u}_i,i=1,\cdots,N\}$ is a Nash equilibrium strategy implies
 $I_1=0 $, 
which is equivalent to
\begin{align*}
&\check{u}_i=-R^{-1}B^T \mathbb{E}[\check{\lambda}_i|\mathcal{F}_t^i]-R^{-1}D^T\mathbb{E}[\check{\beta}^i_i|\mathcal{F}_t^i].
\end{align*}
Thus, we have the optimality system (\ref{eq6a}).
This implies that 
(\ref{eq6a}) admits an adapted solution $(\check{x}_i,\check{\lambda}_i,\check{\beta}_{i}^{j})$.

(ii) If (\ref{eq6a}) admits a solution $(\check{x}_i,\check{\lambda}_i,\check{\beta}_{i}^{j})$, then it can be verified that $I_1=0$,
which with (\ref{eq5a}) implies that $\{\check{u}_i,i=1,\cdots,N\}$ in (\ref{eq7a}) is a set of Nash strategies.
$\hfill \Box$

\emph{Proof of Proposition \ref{prop3.1}.}
When $\Gamma=I$, (\ref{eq9a}) and (\ref{eq9b}) are simplified as
\begin{align}\label{eq12}
\rho K_N=&\dot{K}_N+{A}^TK_N+K_NA-(B^TK_N+D^TK_NC)^T\cr&
\times\Upsilon_N^{-1}(B^TK_N+D^TK_NC)+C^TK_NC,\cr&+\frac{(N-1)^2}{N^2}Q, \quad K_N(T)=0,\\ \label{eq13}
\rho \Pi_N=&
 \dot{\Pi}_N+{A}^T\Pi_N+\Pi_N{A}-\Pi_N B \Upsilon_N^{-1}B^T\Pi_N\cr&-\Pi_N B \Upsilon_N^{-1}(B^TK_N+D^TK_NC)\cr
&-(B^TK_N+D^TK_NC)^T\Upsilon_N^{-1}B^T\Pi_N\cr&-\frac{(N-1)^2}{N^2}Q, \quad  \Pi_N(T)=0.
\end{align}
Note $Q\geq0$. We obtain that $(\ref{eq12})$ admits a solution $K_N\geq 0$ (\cite{SLY16}). Let $P_N=K_N+\Pi_N$. Then the equation
\begin{align}\label{eq14a}
\rho P_N=&\dot{P}_N+{A}^TP_N+P_NA-(B^TP_N+D^TK_NC)^T\Upsilon_N^{-1}\cr&\times(B^TP_N+D^TK_NC)
  +C^TK_NC,\
    P_N(T)=0
\end{align}
admits a unique solution $P_N\geq0$. Since $K_N\geq0$, we obtain that $(\ref{eq13})$ admits a solution, which further implies A2) holds.
$\hfill \Box$

\section{Proofs for Section \ref{sec4}}
\def\theequation{B.\arabic{equation}}
\setcounter{equation}{0}

\emph{Proof of Proposition \ref{lem2c}}. For the case $\Gamma=I$, (\ref{eq15c})-(\ref{eq16c}) can be written as
\begin{align}\label{eq32}
 \rho K_N&=A^TK_N+K_NA-(B^TK_N+D^TK_NC)^T\Upsilon_N^{-1}\cr&\times(B^TK_N+D^TK_NC)+C^TK_NC+\Big(\frac{N-1}{N}\Big)^2 Q,
\\ \label{eq33}
  \rho {P}_N&={A}^TP_N+P_NA-(B^TP_N+D^TK_NC)^T \Upsilon_N^{-1}\cr&\times(B^TP_N+D^TK_NC)+C^TK_NC.
\end{align}
Since A3) holds,
then (\ref{eq32}) admits a stabilizing solution. From \cite{M77}, (\ref{eq33}) admits a stabilizing solution if and only if $ \mathcal{M}_{I}$ has no eigenvalues on the imaginary axis. For the case $C=0$, $ \mathcal{M}_{I}$ has no eigenvalues on the imaginary axis if and only if $A-\frac{\rho}{2}I$ has no eigenvalues in the imaginary axis. Then the theorem follows.
$\hfill \Box$

To prove Theorem \ref{thm11}, we first provide a lemma, which shows uniform stability of the closed-loop systems.
\begin{lemma}\label{lem7}
	Assume that A1), A3), A4) hold and $N$ is sufficiently large such that $I-\frac{1}{N}\Gamma$ is nonsingular.
	Then (\ref{eq17c}) admits a unique solution $s_N\in C_{\rho/2}([0,\infty),\mathbb{R}^n)$
and  
the following holds:
\begin{equation*}
	\sum_{i=1}^N\mathbb{E}\int_0^{\infty} e^{-\rho t} \left(\|\check{x}_i(t)\|^2+\|\check{u}_i(t)\|^2\right)dt<\infty.
	\end{equation*}	
\end{lemma}

\emph{Proof.} 
 In view of A4), $A-B\Upsilon_N^{-1}(B^TP_N+D^TK_NC)-\frac{\rho}{2}I$ is Hurwitz, and hence
$  \int_0^{\infty}e^{-\rho t}(\mathbb{E}x_i)^2dt<\infty.$
From an argument in \cite[Appendix A]{WZ12}, we obtain (\ref{eq17c}) admits a unique solution $s_N\in C_{\rho/2}([0,\infty),\mathbb{R}^n)$.
Denote
$$\bar{A}_N\stackrel{\Delta}{=}A-B\Upsilon_N^{-1}(B^TK_N+D^TK_NC),$$
$$\bar{C}_N\stackrel{\Delta}{=}C-D\Upsilon_N^{-1}(B^TK_N+D^TK_NC),$$
\begin{align*}
	\bar{f}_N\stackrel{\Delta}{=}&f-B\Upsilon_N^{-1}\big[B^T((P_N-K_N){\mathbb{E}}[\check{x}_i]+s_N)+D^TP\sigma\big],\\
	\bar{\sigma}_N\stackrel{\Delta}{=}&\sigma-D\Upsilon_N^{-1}\big[B^T((P_N-K_N){\mathbb{E}}[\check{x}_i]+s_N)+D^TP\sigma\big].
\end{align*}
 After the control (\ref{eq14c}) is applied, we have
\begin{align}\label{eq20i}
d\check{x}_i(t)=&\big[\bar{A}_N\check{x}_i(t)+\bar{f}_N(t)\big]dt+\big[\bar{C}_N\check{x}_i(t)+\bar{\sigma}_N(t) \big]dw_i(t).
\end{align}
Note that $N$ is sufficiently large such that $I-\frac{1}{N}\Gamma$ is nonsingular. By A3) and \cite{ZZC08}, we obtain that $\big[A-\frac{\rho}{2}I, C,\sqrt{\big(I-\frac{1}{N}\Gamma\big)^T Q\big(I-\frac{1}{N}\Gamma\big)}\big]$ is exactly detectable, and hence $(\bar{A}_N-\frac{\rho}{2}I, \bar{C}_N)$ is mean-square stable. Let $Y_N$ satisfy
\begin{align}
	&Y_N(\bar{A}_N-\frac{\rho}{2}I)+(\bar{A}_N-\frac{\rho}{2}I)^TY_N+(\bar{C}_N)^TY_N\bar{C}_N=-2I.
\end{align}
From \cite{ZZC08}, we have $Y_N>0$.
By It\^{o}'s formula and (\ref{eq20i}),
\begin{align}\label{eq40a}
  &\mathbb{E}[e^{-\rho T}\check{x}_i^T(T)Y_N\check{x}_i(T)-\check{x}_i^T(0)Y_N\check{x}_i(0)]\cr
  =&\mathbb{E}\int_0^Te^{-\rho t}\big[\check{x}_i^T(Y_N\bar{A}_N+\bar{A}_N^TY_N+\bar{C}_N^TY_N\bar{C}_N-\rho Y_N)\check{x}_i\cr
  &+2\big(Y_N\bar{f}_N+\bar{C}_N^TY_N\bar{\sigma}_N\big)^T\check{x}_i+\sigma_N^T Y_N\sigma_N\big]dt\cr
  \leq&\mathbb{E}\int_0^Te^{-\rho t}\big[-\check{x}_i^T\check{x}_i+\|Y_N\bar{f}_N+\bar{C}_N^TY_N\bar{\sigma}_N\|^2\cr&+\bar{\sigma}_N^T Y_N\bar{\sigma}_N\big]dt.
\end{align}
From this, 
we have 
\begin{align*}
\mathbb{E}\int_0^{\infty}e^{-\rho t}\|\check{x}_i\|^2dt\leq &\mathbb{E}\int_0^{\infty}e^{-\rho t}\big[\|Y_N\bar{f}_N+\bar{C}_N^TY_N\bar{\sigma}_N\|^2\cr
&+\bar{\sigma}_N^T Y_N\bar{\sigma}_N\big]dt+c_0\leq c.
\end{align*}
This with (\ref{eq14c}) completes the proof. $\hfill \Box$

\emph{Proof of Theorem \ref{thm11}.}
Denote $\tilde{u}_i=u_i-\check{u}_i$ and $\tilde{x}_i=x_i-\check{x}_i$. Then $\tilde{x}_i$ satisfies ($ \tilde{x}_i(0)=0$)
\begin{equation}\label{eq73}
	\begin{split}		 &d\tilde{x}_i=(A\tilde{x}_i+B\tilde{u}_i)dt+(C\tilde{x}_i+D\tilde{u}_i)dw_i.
	\end{split}
\end{equation}
By Lemma \ref{lem7},
\begin{equation}
	\sum_{i=1}^N\mathbb{E}\int_0^{\infty} e^{-\rho t} \left(\|\tilde{x}_i(t)\|^2+\|\tilde{u}_i(t)\|^2\right)dt<\infty.
	\end{equation}	
From (\ref{eq2}), we have $
J_{i}(u_i,\check{u}_{-i})
=J_{i}(\check{u}_i,\check{u}_{-i})+\tilde{J}_{i}(\tilde{u}_i,\check{u}_{-i})+\mathcal{I}_i,
$
where
\begin{align*}
\tilde{J}_{i}(\tilde{u}_i,\check{u}_{-i})\!\stackrel{\Delta}{=}&\mathbb{E}\!\int_0^{\infty}\!\!e^{-\rho t}\!\big[\|\tilde{x}_i(t)-\!\frac{1}{N}\Gamma \tilde{x}_i(t)\|^2_Q+\|\tilde{u}_i(t)\|^2_{R}\big]dt,\\
\mathcal{I}_i=&2\mathbb{E}\int_0^{\infty}e^{-\rho t}\Big[\big(\check{x}_i(t)-\Gamma \check{x}^{(N)}(t)-\eta(t)\big)^T\cr&\times Q\big(\tilde{x}_i(t)-\frac{1}{N}\Gamma \tilde{x}_i(t)\big)+\check{u}_i^T(t)R\tilde{u}_i(t)\Big]dt.
\end{align*}
By applying It\^{o}'s formula with 
(\ref{eq15c}) and (\ref{eq16c}), we have
 \begin{align}\label{eq35}
 0\!=&\limsup_{T\to\infty}\mathbb{E}\big[ e^{-\rho T}\tilde{x}_i^T(T)\big(K_N\check x_i(T)\cr&+(P_N-K_N)\mathbb{E}[\check{x}_i(T)]+s_N(T)\big)\big]\cr =\!&\sum_{i=1}^N\mathbb{E}\!\int_{0}^{\infty}\!\!\!e^{-\rho t} \Big\{\tilde{x}_i^T\big[ (A^TK_N+K_NA-\rho K_N+C^TK_NC\cr&-(BK_N+D^TK_NC)^T\Upsilon_N^{-1}(B^TK_N+D^TK_NC)) \check{x}_i\cr&-\big({A}^T(P_N-K_N)+(P_N-K_N){A}-\rho P_N\cr&-(B^TP_N+D^TK_NC)^T\Upsilon_N^{-1}(B^TP_N+D^TK_NC)\cr&
\!+(B^TK_N\!+D^TK_NC)^T\Upsilon_N^{-1}(B^TK_N\!+D^TK_NC)\big)\mathbb{E}[\check{x}_i]\cr&+\dot{s}_N-\rho s_N+\big[{A}-B\Upsilon_N^{-1}\big( B^T(K_N+\Pi_N)\cr&+D^TK_NC\big)\big]^T s_N+(K_N+\Pi_N)f\cr
&+\big[C-D\Upsilon_N^{-1}\big( B^T(K_N+\Pi_N)+D^TK_NC\big)\big]^TK_N\sigma\big] \cr
&
+\tilde{u}_i^T B^T(K_N\check x_i+(P_N-K_N)\mathbb{E}[\check{x}_i]+s_N)\cr&+\tilde{u}_i^T D^TK_ND\check{u}_i\Big\}dt\cr
= &\mathbb{E}\int_{0}^{\infty}e^{-\rho t}\Big\{- \tilde{x}_i^T\Big[ \big(I-\frac{1}{N}\Gamma\big)^TQ\big(I-\frac{1}{N}\Gamma\big)\check{x}_i\cr&
-\frac{N-1}{N}\big(I-\frac{1}{N}\Gamma\big)Q\Gamma\mathbb{E}[\check{x}_i]-\big(I-\frac{1}{N}\Gamma\big)Q{\eta}\Big]\cr
&-\check{u}_i^TR\tilde{u}_i\Big\}dt.
\end{align}
Note $\check{u}_i,\check{x}_i, \tilde{u}_i,\tilde{x}_i$ are adapted to $\mathcal{F}_t^i$. By the property of conditional expectation,
$$\begin{aligned}
\mathcal{I}_i=&2\mathbb{E}\int_0^{\infty}e^{-\rho t}\mathbb{E}\Big[\big(\check{x}_i(t)-\Gamma
\check{x}^{(N)}(t)-\eta(t)
\big)^TQ\\
&\times\big(\tilde{x}_i(t)-\frac{1}{N}\Gamma
\tilde{x}_i(t)\big)+\check{u}_i^T(t)R\tilde{u}_i(t)\big|\mathcal{F}_t^i\Big]dt\\
=&2\mathbb{E}\int_0^{\infty}e^{-\rho t}\mathbb{E}\Big[\tilde{x}_i^T(t)\big(I-\frac{1}{N}\Gamma\big)^T Q\big(\big(I-\frac{1}{N}\Gamma\big)\check{x}_i(t)\cr&+\check{u}_i^T(t)R\tilde{u}_i(t)-\big(\frac{N-1}{N}\big)\Gamma \mathbb{E}[\check{x}_i(t)]-\eta(t)
\big)^T\Big]dt.
\end{aligned}$$
Comparing this with (\ref{eq35}) leads to $ \mathcal{I}_i=0$. Thus, the theorem follows.  $\hfill \Box$

\emph{Proof of Theorem \ref{thm5.2}.}  Note that (\ref{eq423}) is equivalent to
\begin{align*}
	d(\check{x}_i(t)-\bar{x}_0)=&-\Upsilon_N^{-1}{K}_N(\check{x}_i(t)-\bar{x}_0)dt\cr&-\Upsilon_N^{-1}{K}_N(\check{x}_i(t)-\bar{x}_0)dw_i(t).
\end{align*}
From Definition 
\ref{def1}, we obtain that (\ref{eq423}) reaches mean-square consensus, if $[-B\Upsilon_N^{-1}B^T{K}_N,-D\Upsilon_N^{-1}B^T{K}_N]$ is mean-square stable.
It follows by It\^{o}'s formula that there exists a $c_2>0$ such that
\begin{align*}
&\mathbb{E}[\|\check{x}_i(T)-\bar{x}_0\|^2_{K_N}-\|\check{x}_i(0)-\bar{x}_0\|^2_{K_N}]\cr
=&\mathbb{E}\int_0^T (\check{x}_i(t)-\bar{x}_0)^T\big(-2K_N\Upsilon_N^{-1} K_N\cr&+K_N\Upsilon_N^{-1}K_N\Upsilon_N^{-1} K_N\big)(\check{x}_i(t)-\bar{x}_0)\cr
=&\mathbb{E}\int_0^T (\check{x}_i(t)-\bar{x}_0)^T\big(-K_N\Upsilon_N^{-1}(\Upsilon_N+R)\Upsilon_N^{-1} K_N\big)\cr&\times(\check{x}_i(t)-\bar{x}_0)dt\cr
\leq& -c_2 \mathbb{E}\int_0^T (\check{x}_i(t)-\bar{x}_0)^TK_N(\check{x}_i(t)-\bar{x}_0)dt.
\end{align*}
By Gronwall's inequality, we obtain
$$\mathbb{E}[\|\check{x}_i(T)-\bar{x}_0\|^2_{K_N}]\leq \mathbb{E}[\|\check{x}_i(0)-\bar{x}_0\|^2_{K_N}] e^{-c_2 t},$$
which further gives
$ \mathbb{E}\|\check{x}^{(N)}(t)-\bar{x}_0\|^2\leq c_1 e^{-c_2 t}.$
 \hfill$\Box$

\section{Proofs for Section \ref{sec6}}
\def\theequation{C.\arabic{equation}}
\setcounter{equation}{0}

 \emph{Proof of Theorem \ref{thm1b}.} (Necessity) Suppose $\{\hat{u}_i,i=1,\cdots,N\}$ is a set of social optimal strategies, and $\{\hat{x}_i,i=1,\cdots,N\}$ is the corresponding states of agents.
Let $\{\hat{\lambda}_i, \hat{\beta}_i^j,i,j=1,\cdots,N\}$ be a set of solutions to the second equation of (\ref{eq63}).
 For any $u_i\in \mathcal{U}_{d,i} 
 $ and $\theta\in \mathbb{R}\ (\theta \not= 0)$, let $u_i^{\theta}=\hat{u}_i+\theta v_i$. Denote by $x_i^{\theta}$
the corresponding state under the control $u_i^{\theta}$, $ i=1,2,\cdots,N.$
 Let $y_i=(x_i^{\theta}-\check{x}_i)/\theta$. 
It can be verified that
 $y_i$ satisfies (\ref{eq4aa}).
Then by It\^{o}'s formula, 
for any $i=1,\cdots,N$,
\begin{align}\label{eq69a}
0
=& \sum_{i=1}^N\alpha_i^{(N)}\mathbb{E}\int_0^T e^{-\rho t}  \Big[\big\langle -\big[Q\hat{x}_i-Q_{\Gamma} \hat{x}^{(\alpha)}-\bar{\eta}\big],y_i\big\rangle\cr&+\langle B^T\hat{\lambda}_i+D^T\hat{\beta}_i^i,v_i\rangle\Big] dt.
\end{align}
We have
${J}_{\rm soc,T}(u^{\theta})-{J}_{\rm soc,T}(\hat{u})=2\theta I_1+{\theta^2}I_2,
$
where $u^{\theta}=(u^{\theta}_1,\cdots,u^{\theta}_N)$, $y^{(\alpha)}=\frac{1}{N}\sum_{j=1}^N\alpha_j^{(N)}y_j$  and
\begin{align*}
I_1\stackrel{\Delta}{=}&\sum_{i=1}^N\alpha_i^{(N)}\mathbb{E}\int_0^T  e^{-\rho t}\big[\big\langle Q\big(\hat{x}_i-(\Gamma\hat{x}^{(\alpha)}+\eta)\big),\cr
&y_i-\Gamma y^{(\alpha)}\big\rangle +
\langle R \check{u}_i,v_i\rangle \big]dt,\cr
I_2\stackrel{\Delta}{=}&\sum_{i=1}^N\alpha_i^{(N)}\mathbb{E}\int_0^T e^{-\rho t}\big[\big\|y_i-\Gamma y^{(\alpha)}
   \big\|^2_{Q}
+\|v_i\|^2_{R}\big]dt.
\end{align*}
Note that
{\small \begin{align*}
  &\sum_{i=1}^N\alpha_i^{(N)}\mathbb{E}\int_0^Te^{-\rho t} \big\langle Q\big(\hat{x}_i-(\Gamma\hat{x}^{(\alpha)}+\eta)\big),\Gamma y^{(\alpha)}\big\rangle dt\cr
 =& \sum_{j=1}^N \mathbb{E}\!\int_0^T\!e^{-\rho t}  \Big\langle  {\Gamma^TQ} \sum_{i=1}^N\alpha_i^{(N)}\big(\hat{x}_i-(\Gamma\hat{x}^{(\alpha)}+\eta)\big),\alpha_j^{(N)} y_j\Big\rangle  dt\cr
 =& \sum_{j=1}^N \alpha_j^{(N)}\mathbb{E}\int_0^Te^{-\rho t} \big\langle {\Gamma^TQ} \big((I-\Gamma)\hat{x}^{(\alpha)}-\eta\big), y_j\big\rangle  dt.
\end{align*}}
From this and (\ref{eq69a}), one can obtain that
\begin{align}  \label{eq72}
 I_1
=&\sum_{i=1}^N\alpha_i^{(N)}\mathbb{E}\int_0^Te^{-\rho t} \Big\langle R\hat{u}_i+B^T\hat{\lambda}_i+D^T\hat{\beta}_i^i,v_i\Big\rangle dt\cr
=&\sum_{i=1}^N\alpha_i^{(N)}\mathbb{E}\int_0^Te^{-\rho t} \Big\langle R\hat{u}_i+B^T\mathbb{E}[\hat{\lambda}_i|\mathcal{F}_t^i]\cr&\qquad\qquad+D^T\mathbb{E}[\hat{\beta}^i_i|\mathcal{F}_t^i],v_i\Big\rangle dt.
\end{align}
where the second equation holds
by the fact that $\hat{u}_i,v_i$ are adapted to $\mathcal{F}_t^i$. 
 Since $Q \geq0$ and $R> 0$, we have
 $I_2\geq0$.
From (\ref{eq72}), $\hat{u}=(\hat{u}_1,\cdots,\hat{u}_N)$ is an optimal control of Problem (G$^{\prime}$) if and only if
 $I_1=0 $, 
which is equivalent to
$\check{u}_i=-R^{-1}B^T \mathbb{E}[\hat{\lambda}_i|\mathcal{F}_t^i]-R^{-1}D^T\mathbb{E}[\hat{\beta}^i_i|\mathcal{F}_t^i].$
Thus, we have the following optimality system (\ref{eq63}).
This implies that 
(\ref{eq63}) admits a solution $(\hat{x}_i,\hat{\lambda}_i,\hat{\beta}_{i}^{j})$.

(Sufficiency) On other hand, if (\ref{eq63}) admits a solution $(\hat{x}_i,\hat{\lambda}_i,\hat{\beta}_{i}^{j})$, then it can be verified that $I_1=0$,
which 
implies $\{\hat{u}_i,i=1,\cdots,N\}$ is social optimal control.
$\hfill \Box$

\emph{Proof of Proposition \ref{prop4}.}
Note that
$$Q-\frac{1}{N}Q_{\Gamma}=\frac{N-1}{N}Q+\frac{1}{N}(I-\Gamma)^TQ(I-\Gamma)\geq0.$$ We obtain that $(\ref{eq68})$ admits a solution $\hat{K}_N\geq 0$. Let $\hat{P}_N=\hat{K}_N+\hat{\Pi}_N$. Then
\begin{equation*}
  \begin{aligned}
\rho \hat{P}_N=  &\dot{\hat{P}}_N+{A}^T\hat{P}_N+\hat{P}_NA+C^T\hat{K}_NC-(B^T\hat{P}_N
  \cr&+D^T\hat{K}_NC)^T\hat{\Upsilon}_N^{-1}(B^T\hat{P}_N+D^T\hat{K}_NC)\cr&+(I-\Gamma)^TQ(I-\Gamma),  \hat{P}_N(T)=0
  \end{aligned}
\end{equation*}
admits a unique solution $\hat{P}_N\geq0$. This further implies (\ref{eq69}) that admits a solution.
\hfill $\Box$

{\emph{Proof of Theorem \ref{thm6.2}.}} The proof for social optimality of (\ref{eq75}) is similar to Theorem \ref{thm3}, and so we omit it here. Note that $\hat{x}_i,\hat{u}_i$ are adapted to $\mathcal{F}_t^i$ and $\mathbb{E}\big[\hat{x}_j|\mathcal{F}_t^i\big]=\mathbb{E}\big[\hat{x}_j\big]=
\mathbb{E}\big[\hat{x}_i\big]=\mathbb{E}\big[\hat{x}^{(N)}]$.
By direct calculations,
\begin{align*}
&J_{\rm soc,T}(\hat{u})
\cr
=&  \sum_{i=1}^N\mathbb{E}\int_0^Te^{-\rho t}\mathbb{E}\Big[\big\|\hat{x}_i-\Gamma \hat{x}^{(N)}-\eta\big\|^2_{Q}+\|\hat{u}_i\|^2_{R}\Big|\mathcal{F}_t^i\Big]dt\cr
=&\sum_{i=1}^N\mathbb{E}\int_0^Te^{-\rho t}\Big[\Big\|\big(I-\frac{1}{N}\Gamma\big)\hat{x}_i\Big\|^2_Q\cr&-2x_i^T\big(I-\frac{1}{N}\Gamma\big)^TQ\frac{\Gamma}{N}\sum_{j\not= i}\mathbb{E}[\hat{x}_j]+\frac{1}{N^2}\mathbb{E}\Big\|\Gamma\sum_{j\not= i}\hat{x}_j\Big\|^2_{Q }
\cr&+\|\eta\|_Q^2-2\eta^TQ(I-\Gamma)\hat{x}_i+\|\hat{u}_i\|^2_{R}\Big]dt.
\end{align*}
Note that
\begin{align*}
  &\mathbb{E}\Big\|\Gamma\sum_{j\not= i}\hat{x}_j\Big\|^2_{Q }=\sum_{j\not= i}\sum_{k\not= i} \mathbb{E}[\hat{x}_j^T\Gamma^T{Q }\Gamma\hat{x}_k]\cr
  =&\sum_{j\not= i}\mathbb{E}[\hat{x}_j^T\Gamma^T{Q }\Gamma\hat{x}_j]+\sum_{j\not= i}\sum_{k\not=j, i}  \mathbb{E}[\hat{x}_j]^T\Gamma^T{Q }\Gamma\mathbb{E}[\hat{x}_k]\cr
  =&(N-1)\mathbb{E}[\hat{x}_i^T\Gamma^T{Q }\Gamma\hat{x}_i]\cr
  &+(N-1)(N-2)  \mathbb{E}[\hat{x}_i]^T\Gamma^T{Q }\Gamma\mathbb{E}[\hat{x}_i].
\end{align*}
Then we further have
\begin{align*}
&J_{\rm soc,T}(\hat{u})
\cr=&\sum_{i=1}^N\mathbb{E}\int_0^Te^{-\rho t}\Big[\Big\|\big(I-\frac{1}{N}\Gamma\big)\hat{x}_i\Big\|^2_Q+\frac{N-1 }{N^2}\mathbb{E}[\hat{x}_i^T\Gamma^T{Q }\Gamma\hat{x}_i]\cr&-2\frac{N-1}{N}x_i^T\big(I-\frac{1}{N}\Gamma\big)^TQ\Gamma\mathbb{E}[\hat{x}_i]
\cr&+\frac{(N-1)(N-2)}{N^2}  \mathbb{E}[\hat{x}_i]^T\Gamma^T{Q }\Gamma\mathbb{E}[\hat{x}_i]+\|\eta\|_Q^2\cr&-2\eta^TQ(I-\Gamma)\hat{x}_i+\|\hat{u}_i\|^2_{R}\Big]dt\cr
=&\sum_{i=1}^N\mathbb{E}\int_0^Te^{-\rho t}\Big(\|\hat{x}_i-\mathbb{E}[\hat{x}_i]\|^2_{Q-Q_{\Gamma}/N}+\|\mathbb{E}[\hat{x}_i]\|^2_{Q-Q_{\Gamma} }\cr&+\|\eta\|_Q^2-2\eta^TQ(I-\Gamma)\mathbb{E}[\hat{x}_i]+\|\hat{u}_i-\mathbb{E}[\hat{u}_i]\|^2_{R}\cr&+\|\mathbb{E}[\hat{u}_i]\|^2_{R}\Big)dt\cr
=&\sum_{i=1}^N\mathbb{E}\big[\|x_{i0}\!-\mathbb{E}[\hat{x}_i(0)]\|^2_{\hat{K}_N}\!+
\|\mathbb{E}[\hat{x}_i(0)]\|^2_{\hat{P}_N}\cr&+2s^T(0)\mathbb{E}[\hat{x}_i(0)]\big]
+\!\sum_{i=1}^N\mathbb{E}\int_0^Te^{-\rho t}\Big(\big\|\hat{u}_i\!-\mathbb{E}[\hat{u}_i]\cr&+\hat{\Upsilon}_N^{-1}B^T\hat{K}_N(\hat{x}_i\!-\mathbb{E}[\hat{x}_i])\big\|^2_{\hat{\Upsilon}_N}\cr
&+\big\|\mathbb{E}[\hat{u}_i]+\hat{\Upsilon}_N^{-1}B^T\hat{P}_N \mathbb{E}[\hat{x}_i]\big\|^2_{\hat{\Upsilon}_N}\Big)dt+  Nq_T^N\cr
= &\sum_{i=1}^N\mathbb{E}\big[\|x_{i0}-\mathbb{E}[{x}_i(0)]\|^2_{P_N}+\|\mathbb{E}[{x}_i(0)]\|^2_{\Pi_N}\cr&+2s^T_N(0)\mathbb{E}[{x}_i(0)]\big]+ Nq_T^N,
\end{align*}
where $q_T^N$ is given by (\ref{eq17a}).    \hfill{$\Box$}

\end{document}